\documentclass{marticle}
\begin{document}
\title{Rotating Walks, Derandomized Squares, and Permutation Conjugated Walks}
\author{
    Vedat Levi
    Alev\thanks{\href{mailto:vedata@ariel.ac.il}{vedata@ariel.ac.il},
        Ariel University.
    }
    \and
    Shravas Rao\thanks{
        \href{mailto:shravas@pdx.edu}{shravas@pdx.edu},
        Portland State University.
    }
}
\date{\today}
\maketitle
\begin{abstract}
    Using the recent strong convergence results of Bordenave and
    Collins (Annals of Mathematics) and Chen et
    al.~(Annals of Mathematics), we study the average case behavior
    of rotating random walks of Cohen and Maor (STOC 2023) and the
    derandomized square construction of Rozenman and Vadhan (RANDOM 2005).

    Concretely, we prove a conjecture of Cohen et al.~(ITCS 2025)
    regarding the average case behavior of derandomized squaring for
    random graphs and show that the existential bounds of Cohen and
    Maor for their rotating random walk construction hold with
    high-probability for random graphs. Using the results of
    O'Donnell and Wu (FOCS 2020), we can in polynomial time construct
    a family of graphs and rotations which almost saturate the bounds
    of Cohen and Maor, and meet the bounds conjectured by Cohen et
    al.~for all vertex-transitive auxiliary graphs of fixed degree in
    the derandomized square.

    Finally, we study independent random permutation conjugates of
    deterministic regular graph adjacency matrices. We prove
    quantitative weak convergence to freely independent tree
    adjacency operators when the degrees are fixed and the girths
    tend to infinity. Separately, a covering argument gives
    deterministic spectral lower bounds for polynomials with
    non-negative coefficients,
    recovering the lower bound of Cohen et al.~for rotating walks.
\end{abstract}
\section{Introduction}\label{sec:introduction} 
\paragraph{Spectral Graph Theory Beyond Worst Case Analyses.} Recall
that an $(N, d, \lambda)$-graph $G = (V, E)$ is a $d$-regular graph
on $N := |V|$ vertices such that $\lambda(G) \le \lambda$, i.e., the
largest magnitude of any non-trivial eigenvalue of the adjacency
matrix $A$ of $G$ is
bounded by $\lambda$. Expander graphs are highly desirable objects
because they maintain high connectivity even when $d = O(1)$, because
of this they have found many applications in theoretical computer
science. The general idea is a random walk over an expander graph
following the power $G^t$ can be thought as being very similar to
independent sampling. We refer to \cite{HooryLW06,Vadhan12} for a
detailed account of the
applications of expander graphs.

The celebrated Alon-Bopanna bound
\cite{Nilli91} states that we necessarily have $\lambda \ge
2\sqrt{d-1} - o(1)$ and graphs attaining $\lambda = 2\sqrt{d-1}$ are
called Ramanujan graphs.
Although Ramanujan graphs are abundant \cite{HuangMY24} and explicit
constructions exist algebraic constructions exist
\cite{LubotzkyPS86,Morgenstern94}, simple operations
such as powering often disturb the Ramanujan quality that a graph
has. Recalling that the adjacency matrix of the graph $G^t$ is $A^t$,
we have that
\[ \textrm{If}~G~\textrm{is an}~(N, d,
    2\sqrt{d-1})\textrm{-graph then}~G^t~\textrm{is an}~(N, d^t, 2^t
(d-1)^{t/2})\textrm{-graph}\]
where we see that the eigenvalue $\lambda^t = 2^t (d-1)^{t/2}$ of
$G^t$ is almost a factor of $2^t$ away from being Ramanujan. This
identifies a bottle-neck for algorithms using random walks.

In \cite{CohenM23} a solution is proposed to overcome this
bottle-neck. Let $\Pi_1,\ldots, \Pi_{k-1}$ be permutation matrices.
The rotating walk over $G$ with the permutation matrices $\rv \Pi = (\Pi_1,
\ldots, \Pi_{k-1})$ is the which takes $k$ steps in the graph $G$
which are interleaved by permutation steps according to
$\Pi_1,\ldots, \Pi_{k-1}$. Though their result holds for general
high-girth graphs, for simplicity we focus on the case where $G$ is a
high-girth Ramanujan graph. For technical reasons, they define
\begin{equation}\label{eq:rotwalk_show}
    \Arot(\rv \Pi) = (A \cdot
        \Pi_1 \cdot A \cdot \Pi_2 \cdot A \cdots
    A \cdot \Pi_{k-1} \cdot A)(A \cdot \Pi_1 \cdot A \cdot \Pi_2 \cdot A \cdots
    A \cdot \Pi_{k-1} \cdot A)^*.
\end{equation}
and shows that there exists a sequence of permutation matrices such
that
\begin{equation}\lambda\parens*{\Arot(\rv \Pi)} \le \parens*{\frac{k+1}{k}}^k
    \cdot (k+1) \cdot d^k + o(1)\label{eq:rotbd-show}
\end{equation}
when $k$ is small relative to the girth of $G$, with a similar result
extending a similar conclusion to non-Ramanujan high-girth graphs.
The main technique used in \cite{CohenM23} is \emph{finite-free
probability} \cite{MarcusSS22,Marcus21}. In \cite{CohenM23} a
suspicion is voiced that their
worst case results capture the average case of the rotating walk.

The result \eqref{eq:rotbd-show} is shown to be tight in a follow-up
paper by \cite{CohenCMP25}, by showing that for any collection of
permutation matrices $\rv P = (P_1,\ldots, P_k)$ we necessarily have
\begin{equation}\label{eq:rot_lbd} \lambda\parens*{\Arot(\rv P)} \ge
    \parens*{\frac{k+1}{k}}^k
    \cdot (k+1) \cdot (d-1)d^{k-1} - o(1).
\end{equation}
See a more detailed discussion in
\cref{sec:rotating_infinite_trees_and_expanders}.

\cite{CohenCMP25} also analyze the derandomized square construction
of \cite{RozenmanV05}. This construction was introduced by
Rozenman and Vadhan to simplify the analysis of the seminal result of
Reingold \cite{Reingold08}, which using the zig-zag product
construction \cite{ReingoldVW02} showed that undirected connectivity
problem can be decided in logarithmic space. In intuitive terms,
given an $(N, d)$-graph $G$ and an $(d, k)$-graph $H$, the
derandomized square $S = G \dsq H$ is obtained by imposing a copy of
$H$ in the neighborhood of every vertex $v \in V$. Writing
$\omega(\bullet)$ for the normalized expansion of a graph, i.e.,
$\omega(X) = \lambda(X)/D$ for any $D$-regular graph $X$, they prove
\begin{equation}
    \label{eq:dsq_ub_show}
    \omega(G \dsq H) \le 1 - (1- \omega(G)^2)(1- \omega(H)) \le
    \omega(G)^2 + \omega(H).
\end{equation}
In contrast \cite{CohenCMP25} using techniques of analytic
combinatorics (in the sense of \cite{FlajoletS09}) proves that
\begin{equation}
    \lambda(G \dsq H) \ge \Lambda_H - o_N(1)\qquad\textrm{where}\qquad
    \Lambda_H =
    d \cdot \parens*{ x_0 - (d-1) \frac{\chi_{x_0}(H)}{\chi'_{x_0}(H)} }
    \label{eq:lowerbound_dsq_show}
\end{equation}
where $\chi_\bullet(H)$ denotes the characteristic polynomial of the
adjacency matrix $A_H$ of $H$ evaluated at $x_0$, which is the
largest root to the equation
\begin{equation}
    (d-1) \cdot \chi_x(H) \cdot \chi''_x(H) = (d-2) \cdot \chi'_x(H)^2.
    \label{eq:x0_def_show}
\end{equation}
They conjecture that the bound \eqref{eq:lowerbound_dsq_show} is
saturated by random graphs from the matching model, and prove that
there are $H$-local graphs (graphs where every vertex neighborhood
contains a copy of $H$) saturating \eqref{eq:lowerbound_dsq_show}
(\cite[Conjecture 2.6]{CohenCMP25}).
Crucially, while every derandomized square is $H$-local not all
$H$-local graphs are derandomized squares.
See a more detailed discussion in
\cref{sec:average_behavior_of_derandomized_squaring}.

Our main insight in this work is that both of the aforementioned
results on rotating walks and derandomized square can be analyzed from
the point of view strong convergence for random permutation and
matching matrices, which permits a deeper understanding and
quantitative strengthening of the results.

\paragraph{Strong Convergence and Free Probability Theory.}
To analyze the rotating walk and the derandomized square
constructions more carefully, we will borrow tools from free
probability \cite{NicaS06, Speicher25} and strong convergence of
random permutation and matching matrices, \cite{BordenaveC19,
OdonnellW20, ChenGVTvH26}. For completeness, we give a brief summary.

Free probability allows studying spectral properties of random
matrices by studying the spectral properties of free
objects which limit them. Voiculescu's seminal work
\cite{Voiculescu91} showed that (mixed) moments
of iid random GUE matrices converge to the moments of free
semi-circular elements. A subsequent result of
\citeauthor{HaagerupT05} \cite{HaagerupT05} established a strong
convergence result showing that the operator norm of arbitrary
non-commutative polynomials in random GUE matrices converge in
probability to the
operator norm of the corresponding non-commutative polynomial in free
semi-circulars. This was extended by Male \cite{Male12} for
combinations of Gaussian matrices and deterministic matrices
(admitting a strong limit) and
\citeauthor{CollinsM14} \cite{CollinsM14} to combinations of iid Haar-random
unitary matrices and deterministic matrices (admitting a strong limit).

On the combinatorial side, the work of Nica \cite{Nica93} established
that mixed moments of iid permutation matrices converge in
probability to the moments of
freely independent Haar unitaries, which was extended by Neagu \cite{Neagu07} to
permutations with restricted cycle length for an appropriate free
family. \citeauthor{BordenaveC19}
\cite{BordenaveC19} then established a strong convergence result that
showed that the (after projecting away from constant functions) norms
of non-commutative polynomial in independent
permutation and perfect matching matrices also converge in probability to
the norms of the corresponding polynomials in free objects.
\citeauthor{OdonnellW20}
\cite{OdonnellW20} showed that
the results of \cite{BordenaveC19} can be constructivized for fixed
parameters, i.e., when
we have a class of low-degree polynomials $\mathcal P$ with bounded
coefficients, we can efficiently find a mixture of permutation and
matching matrices which have a similar spectrum to the corresponding
free objects for all the polynomials in $\mathcal P$ (excluding the
eigenvalue corresponding to the constant function). Finally,
\citeauthor{ChenGVTvH26} \cite{ChenGVTvH26} gave a non-asymptotic
version of strong
convergence for random permutation matrices by using the
\emph{polynomial method}. We refer to the surveys
\cite{vanHandel26,vanHandel26a} for more information on these topics
and \cref{sec:preliminaries}.

\paragraph{Our Results.}
We
resolve the conjecture by \cite{CohenCMP25} which stated that
the inequality \eqref{eq:lowerbound_dsq_show} is saturated by random
graphs from the permutation model. Concretely, we show
\begin{theorem}[Restatement of \cref{thm:dsq_conj}]\label{thm:dsq_conj_show}
    For every $k$-regular vertex-transitive graph $H = ([d], E)$
    (where $d \ge 3, k \ge 1)$ and
    even $N \in \NN$, we have
    \[ \Pp_{G}\sqbr*{ \lambda\parens*{ G
    \dsq H } \ge \Lambda_H + \ee} = o_{N, d}(1) \]
    where $G$ is a random $d$-regular graph drawn from the matching
    model, i.e.,  the adjacency matrix of $G$ satisfies $A = \sum_{i
    = 1}^d M_i$ for iid uniform perfect matching matrices
    $M_1,\ldots, M_d \in \RR^{N \times N}$.
    and $\Lambda_H$ is as
    defined in \eqref{eq:lowerbound_dsq_show}.
\end{theorem}
To establish this result, we first closely trace the proof of
\eqref{eq:dsq_ub_show} by \citeauthor{RozenmanV05}
\cite{RozenmanV05}. This will establish that the derandomized square
$G \dsq H$ can be thought as a quadratic non-commutative polynomial
in perfect matching matrices. Then, we will give an analysis of the
limiting object for the derandomized square. This object can be
identified with the derandomized square $\Tre_d \dsq H$. By adapting
the analysis of \cite{CohenCMP25} for the pruned version $\mathbb X$
of this graph considered in their paper, we will establish a norm
bound for $\Tre_d \dsq H$ and the conjecture will follow.

The strength of this method is that in conjunction with
\cite{OdonnellW20} it will allow us to assert the existence of an
infinite family of base
graphs $G_{N'}$ which are approximately Ramanujan and for each $H$
they almost saturate the lower bound \eqref{eq:lowerbound_dsq_show}. Concretely,
\begin{corollary}[Restatement of \cref{cor:dsq_constr}]
    For any given $N, d, k \in \NN$ ($d \ge 3, k \ge 1$) and $\ee >
    0$, there is a $d$-regular graph $G_{N'} = (V_{N'},
    E_{N'})$ on $N' \in [N, (1+ o_N(1))N]$ vertices such that
    \begin{enumerate}[(i)]
        \item $\lambda(G_{N'}) \le 2\sqrt{d-1} + \ee$, i.e., $G_{N'}$
            is almost Ramanujan.
        \item $\lambda(G_{N'} \dsq H) \le \Lambda_H + \ee$
            simultaneously for every $k$-regular vertex transitive
            graph $H = ([d], E')$, where $\Lambda_H$ is as defined in
            \eqref{eq:lowerbound_dsq_show}.
    \end{enumerate}
    Moreover, the graph $G_{N'}$ can be constructed in time $N^{O_{d, \ee}(1)}$.
\end{corollary}
This result stands in contrast with the $H$-local graphs constructed
by \cite{CohenCMP25} which saturate the bound
\eqref{eq:lowerbound_dsq_show}, which need not always be a
derandomized square.

Our next result concerns the rotating walk from
\eqref{eq:rot_walk_def}.
\begin{theorem}[Simplified Version of
    \cref{thm:permutation_model_nonasymp}]\label{thm:permutation_show}
    Let $G$ be an $(N, d)$-graph drawn from either the permutation
    model (assume in this case $d = 2r$ is even) or the matching model
    (assume in this case that $N$ is even), i.e., suppose that the
    adjacency matrix $A \in \RR^{N \times N}$ of $G$ satisfies
    \[ A = \sum_{i = 1}^r (P_i + P_i^*)\qquad\textrm{or}\qquad A =
    \sum_{i = 1}^d M_i\]
    where $P_i, M_j$ are independent uniform permutation or matching
    matrices for $i \in [r], j \in [d]$. Then, writing $\Arot(\rv
    \Pi)$ for the rotating walk
    over $G$ as defined in \eqref{eq:rotwalk_show}
    \[ \Pp_{G, \rv\Pi}\sqbr*{ \lambda(\Arot(\rv \Pi)) \ge
            \parens*{\frac{k+1}{k}}^k
    \cdot (k+1) \cdot (d-1) d^{k-1} + \ee} = o_{N, d, k}(1).\]
    In particular, when $G$ is drawn from the permutation model, we
    have the following sharper estimate
    \[ \Pp_{G, \rv\Pi}\sqbr*{ \lambda(\Arot(\rv \Pi)) \ge
            \parens*{\frac{k+1}{k}}^k
        \cdot (k+1) \cdot (d-1) d^{k-1} + \ee} \lesssim \frac{1}{N} \cdot
    \parens*{\frac{d^{2k} 4k \log(r+k)}\ee}^8 \cdot \log \frac{e d^{2k}}\ee.\]
\end{theorem}
This result is motivated by the realization that when $G$ is a random
graph drawn from the permutation or matching models, we can think of
the adjacency matrix of the rotating walk as a non-commutative
polynomial in permuation and matching matrices. Thus, by using the
aforementioned results of \cite{BordenaveC19} and \cite{ChenGVTvH26},
this results will follow by a bound on the norm of the limiting
object.

It turns out the limiting object for the rotating walk over a random
graph has a fairly interesting description, namely, it is the
rotating random walk over the infinite $d$-regular tree. The main
tool we use to bound the operator norm of this infinite rotating tree
is the $S$-transform of Voiculescu \cite{Voiculescu87}. This is in
contrast to the proof of \eqref{eq:rotbd-show} in \cite{CohenM23}, which
uses methods of finite free-probability as in \cite{MarcusSS22,
Marcus21}. We note that using the $S$-transform to determine spectral
edges and norms of free operators has some precedents, see
e.g.~\cite{Kargin08,MlotkowskiNPK15,PenningtonSG17}.

The method of \cite{CohenM23} allows one to assert the existence of a
tuple of permutation matrices $\rv \Pi$ for which
\eqref{eq:rotbd-show} holds regardless of the choice of the base
graph, so long as it satisfies a girth and expansion criterion,
nevertheless it does not tell us how to construct such matrices. Our
result works for random graphs, but nevertheless using the methods of
\cite{OdonnellW20} for bounded $k$ and $d$ we can efficiently
construct a family of approximately Ramanujan base graphs $G_{N'}$
for infinitely many $N'$
and an accompanying sequence of permutations $\rv \Pi^{N'} =
(\Pi_1^{N'}, \ldots, \Pi^{N'}_{k-1})$ such that \eqref{eq:rotbd-show} is
approximately satisfied. Namely,
\begin{corollary}[Restatement of \cref{cor:rot_explicit}]
    For each $N, d, k \in \NN$, and $\ee > 0$ there exists some $N'
    \in [N, (1+ o_N(1))
    N]$ such that there exists a $d$-regular graph $G_{N'} =
    (V_{N'}, E_{N'})$ with adjacency matrix $A^{N'}$ and a sequence
    of permutation matrices $\rv
    \Pi^{N'} = (\Pi^{N'}_0,\ldots, \Pi^{N'}_{k-1})$ such that
    \begin{enumerate}[(i)]
        \item $\lambda(A_{N'}) \le 2\sqrt{d- 1} + \ee$,
        \item $\lambda\parens*{\Arot(\rv \Pi^{N'})} \le
            \parens*{\frac{k+1}k}^k (k+1) (d-1)d^{k-1} + \ee.$
    \end{enumerate}
    Moreover, there is a deterministic algorithm which can compute
    the graph $G_{N'}$ and the permutations $\rv \Pi^{N'}$ in time
    $N^{O_{d, k, \ee}(1)}$.
\end{corollary}
We note that a construction with weaker quantitative guarantees already
occurred in the work of \citeauthor{HsiehMZ26} \cite[Theorem
3.6]{HsiehMZ26}. They construct
a $d$-regular base graph $G$ with a sequence of permutations which
satisfy $\lambda\parens*{\Arot(\rv \Pi)} \le O(k^2) \cdot d^k + \ee$
by appealing to Haagerup's inequality (rapid decay property)
\cite{Haagerup79}. By the same strong convergence arguments we use this can be
turned into a high probability norm bound around the weaker margin
$O(k^2) d^k$, whereas both our concentration and construction results
work around the sharp threshold $(k+1)^{k+1}/k^k (d-1)d^{k-1} = O(k) d^k$.

Finally, in order to understand the rotating walk construction for
deterministic base graphs, we study a model of random graphs/matrices which
we call \emph{permutation conjugated random matrices}. Here, we take
a sequence of $d_i$-regular graphs on $N$-vertices $G_{i}$ with
adjacency matrices $A_i$ for all $i = 1,\ldots, r$ and define matrices
\[ B_i = P_i A_i P_i^*\qquad\forall~i \in [r].\]
where $\rv P =(P_1,\ldots, P_r)$ is a sequence of iid uniform
permutation matrices. Intuitively, each $B_i$ is the adjacency matrix
of a graph isomorphic to $G_i$ whose vertex labels are relabeled by a
random permutation. We explain the connection between this model and
the rotating walk construction in
\cref{sec:permutation_conjugated_adjacency_matrices} and give the
formal correspondence in \cref{prop:rotating_redux}.

Although we were unable to show strong convergence results for this
model, we will show a weak convergence result when the graphs have
increasing girth.
\begin{theorem}[Informal Restatement of
    \cref{thm:weak_convergence}]\label{thm:weak_convergence_show}
    Suppose we are given a sequence of $d_i$-regular graphs $G_i$
    satisfying $\girth(G_i) = \omega(1)$ for each $i = 1,\ldots, r$.
    Then, the sequence of matrices $\rv B = (B_1,\ldots, B_r)$ weakly
    converge to $\rv t = (t_1,\ldots, t_r)$ where each $t_i$ is a
    copy of the infinite $d_i$-regular tree $\Tre_d$.

    Equivalently, we have for any fixed $\ee > 0$, and
    non-commutative $*$-polynomial $P$
    \[ \Pp\sqbr*{ \Abs*{\tr P(\rv B) - \tau(P(\rv t))} \ge \ee
    } = o_N(1)\]
    where $\tr(\bullet)$ denotes the normalized trace of $N \times N$
    matrices $\tau(\bullet)$ the trace function of the reduced algebra in
    which $\rv t$ are defined.
\end{theorem}
A simple consequence of \cref{thm:weak_convergence_show} is this
following lower bound,
\begin{corollary}[Informal Restatement of \cref{cor:weak_bd}]
    Suppose we are given a sequence of $d_i$-regular graphs $G_i$
    satisfying $\girth(G_i) = \omega(1)$ for each $i = 1,\ldots, r$.
    Then, the sequence of matrices $\rv B = (B_1,\ldots, B_r)$ weakly
    converge to $\rv t = (t_1,\ldots, t_r)$ where each $t_i$ is a
    copy of the infinite $d_i$-regular tree $\Tre_d$.

    Then,
    \[ \norm*{\left. P(\rv B)\right|_{\perp 1}} \ge \norm*{P(\rv t)} -
    o_N(1)\qquad\textrm{with probability}\qquad 1- o(1).\]
\end{corollary}
When $P$ has only real, non-negative coefficients, we can get a lower
bound which works with probability 1
\begin{theorem}[Informal Restatement of
    \cref{thm:general_lowerbound}]\label{thm:gen_lbd}
    Let $Q$ be fixed a
    non-commutative $*$-polynomial with non-negative real
    coefficients. Then, for every $N \ge 2$, we have
    \[ \norm*{Q(\left. \rv B^N\right|_{\perp 1})} \ge \norm*{Q(\rv t)} -
    o_N(1) \]
    and we assume $d_i=  O(1)$ for all $i = 1,\ldots, r$ with probability 1.
\end{theorem}
In conjunction with the argument that rotating walks can be cast in
the permutation conjugated adjacency matrices model
(\cref{prop:rotating_redux}), this result yields an alternative proof
of \cref{eq:rot_lbd}.

\paragraph{Other Related Work.} Motivated by similar concerns as
\cite{CohenM23}, the work of \cite{HsiehMZ26} uses tools
from free-probability and the results of \cite{OdonnellW20} to design
explicit $\ee$-balanced codes close to the Gilbert-Varshamov bound.
They use walks over a collection of auxiliary graphs obtained by the
work of \citeauthor{OdonnellW20} to show that one can get close to
the Gilbert-Varshamov bound with a simple construction.

The work of \cite{CohenCM24} studies the average case behavior of
the zig-zag product construction of \cite{ReingoldVW02} in the
configuration model.

The work of \citeauthor{WangLZ26} \cite{WangLZ26} is similar in
spirit to \cite{OdonnellW20}. Whereas the results of \cite{OdonnellW20} is an
effective constructivization of the results of \cite{BordenaveC19},
\cite{WangLZ26} considers the task of constructivizing the matrix
concentration results for \cite{BandeiraBMvH23,BrailovskayavH26} for
sums of independent random matrices. These results are effective in a
regime complementary to what we consider, where the graph degree
grows rapidly in graph size.
\paragraph{Acknowledgements.} We would like to thank Gal Maor and
Alon Dogon for discussions regarding their independent work on a
similar topic; and Fernando Granha Jeronimo for his participation in
the earlier stages of this project.
\subsection*{AI Disclosure.} We have started this project by
noticing that the strong convergence techniques of \cite{ChenGVTvH26}
could be used for giving non-asymptotic tail bounds for the quality
of expansion in rotating expanders for random graph models. In
investigating the lower bound
for the rotating walk construction of \cite{CohenM23} by
\cite{CohenCMP25}, we came across the lower bound for the
derandomized square \eqref{eq:lowerbound_dsq_show} and the conjecture
that it should be saturated by random graphs from the matching model.
Realizing that the derandomized square construction can be described
as a non-commutative polynomial following the proof of
\cite{RozenmanV05}, we noticed that strong convergence techniques can
be used to prove this conjecture also.

We have received help from LLM models (ChatGPT-6.1 Sol, ChatGPT-6
Astra) in formalizing some of our intuition and to a very limited capacity
in the preparation of this manuscript.

In particular, some details of the
$S$-transform method which led to the proof of \cref{thm:permutation_show}
was suggested to us by LLMs, when we prompted it for help with a
free-probability proof of the bound \eqref{eq:rotbd-show}. This
materially affected the results of
\cref{sec:rotating_infinite_trees_and_expanders}. Similarly,
we have received help from LLM models, in trying to simplify the
proof of the operator norm bound for which lead to the derandomized
square bound \cref{thm:permutation_show}. This materially affected
the results in \cref{sec:average_behavior_of_derandomized_squaring}.
Finally, we have used LLMs to fix bugs in a previous proof of the
weak convergence result (\cref{thm:weak_convergence}) in
\cref{sec:permutation_conjugated_adjacency_matrices}.

We verified the correctness and originality of all content, including
references.
\section{Preliminaries}\label{sec:preliminaries} 
\subsection{Free Probability Basics}
\label{sub:free_probability_basics} 
We will recall some basic definitions and facts from free probability
theory following \cite{NicaS06}.

\begin{defn}[$C^*$-probability spaces]
    A non-commutative probability space is a pair $(\Adj, \tau)$
    where $\Adj$ is a unital algebra over $\CC$ and a linear
    functional $\tau: \Adj \to \CC$ such that $\tau(1_{\Adj}) = 1$.
    If $\tau(ab) = \tau(ba)$ for all $a, b \in \Adj$ we will call
    $\tau$ a trace and $(\Adj, \tau)$ tracial.

    $\Adj$ is called a $*$-algebra if there exists an
    antilinear $*$-operation over $\Adj$ such that $\Adj \ni a
    \mapsto a^* \in \Adj$, such that for all $a, b \in \Adj$ we have
    $(a^*)^* = a$ and $(ab)^* = b^* a^*$. If $\tau(a^* a) \ge 0$, we will
    call the trace $\tau$ positive and the pair $(\Adj, \tau)$ a
    $*$-probability space. When $\tau(a^* a) = 0$ implies that $a =
    0$, we will call $\tau$ faithful.

    $\Adj$ is called a $C^*$-algebra if it is equipped with a
    complete norm $\norm*{\bullet}$, i.e., a norm satisfying
    $\norm*{a^* a} = \norm*{a}^2$ and $\norm*{a b} \le \norm*{a}
    \norm*{b}$ for all $a, b \in \Adj$.

    Finally, a $C^*$-probability space $(\Adj, \tau)$ is a pair where
    $\Adj$ is a unital $C^*$-algebra $\Adj$ for which $\tau$ is a
    positive trace.
\end{defn}

We note that $C^*$-algebras $(\Adj, \tau)$ can be equipped with a
norm $\norm*{\bullet}$ such
that $\norm*{a^* a} = \norm{a}^2$ for any $a \in \Adj$ and $\norm*{a
b} \le \norm*{a} \norm*{b}$ for all $a, b\in\Adj$.
In particular, this norm satisfies the
following formula
\begin{lemma}[\cite{NicaS06}, Proposition 3.17]\label{lem:norm_def}
    Let $(\Adj, \tau)$ be a $C^*$-probability space such that $\tau$
    is faithful. Then,
    \[ \norm*{a} = \lim_{p \to \infty} \parens*{\tau\parens*{(a^*
    a)^k}}^{\frac{1}{2k}}\]
\end{lemma}
\begin{xmpl}[Matrix Algebra]\label{xmpl:matrix}
    Let $\Mat_N(\CC)$ denote the collection of $N \times N$ matrices
    over $\CC$, where $A^*$ denotes the Hermitian conjugate of $A$,
    i.e.~$A^* = \bar{A}^\top$ and $\tr(A) := \tr_N(A) = N^{-1}
    \sum_{i = 1}^N A_{ii}$ the normalized trace. Then, $(\Mat_N(\CC),
    \tr)$ is a $C^*$-probability space, and the corresponding norm
    $\norm{\bullet}$ is simply the operator norm.
\end{xmpl}
\begin{xmpl}[Reduced $C^*$ Algebras]\label{xmpl:group_algebras}
    Let $G$ be a discrete group with identity element $e \in G$, we
    will write $\ell_2(G)$ for the vector space
    \[ \ell_2(G) = \set*{ \xi : G \to \CC \col\sum_{g \in G}
    |\xi(g)|^2 < \infty }.\]
    We will write $B(\ell_2(G))$ for the collection of bounded
    operators over $\ell_2(G)$, and recall that the left-regular
    representation $\lambda : G \to B(\ell_2(G))$ is defined via the relation
    \[ \lambda(g) \delta_h = \delta_{gh}\qquad\forall~g, h \in G\]
    where $\delta_g$ denotes the indicator function of $g$,
    i.e.~$\delta_g(h) = \One[g = h]$ for all $g, h \in G$.

    We will write, $C^*_{\mathrm{red}}(G)$ for the space
    \[ C^*_{\mathrm{red}}(G) =
    \mathrm{cl}\parens*{\sspan\set*{\lambda(g) \col g \in G}}\]
    where $\mathrm{cl}(\bullet)$ denotes the closure in the norm
    topology. We will equip this space
    with a trace $\tau := \tau_e$ by setting
    \[ \tau\parens*{a} = \inpr*{a \delta_e, \delta_e} = \sum_g
    [a\delta_e](g) \cdot \delta_e(g)\qquad\forall~a \in C^*_{\mathrm{red}}(G).\]
    It is easy to verify that this renders $C^*_{\mathrm{red}}(G)$
    into a $C^*$-probability space. Henceforth, we will call
    $C^*_{\mathrm{red}}$ the reduced $C^*$-algebra of the group $G$.
\end{xmpl}

\begin{defn}[spectrum, spectral distribution]\label{defn:spec_def}
    The spectrum $\Spec(a)$ of an element $a \in \Adj$ is
    defined to be,
    \[ \Spec(a) = \Big\{z\in \CC \col z \cdot 1_{\Adj} -
    a~\textrm{is not invertible}\Big\}.\]

    When $a \in \Adj$ is a normal element, i.e., $a^* a = a a^*$ we
    define the spectral measure $\mu_a$ of $a$ by the unique
    probability measure satisfying
    \[ \int_{\CC} z^p \bar z^q~\dif\mu_a(z) = \tau\parens*{a^p
    (a^*)^q}\qquad\forall~p, q \in \NN.\]
\end{defn}

We can verify easily that when $(\Adj, \tau)$ is the matrix algebra
from \cref{xmpl:matrix}, the spectrum $\Spec(M)$ of any matrix $M \in
\CC^{N \times N}$ simply coincides with the spectrum of $M$ from
linear algebra, i.e., the collection of eigenvalues of $M$.
Similarly, when $M$ is normal, the spectral distribution $\mu_M$ is given by the
empirical eigenvalue distribution $\mu_M(x) = \frac{1}{N}
\sum_{\lambda \in \Spec(M)} \One[x = \lambda]$ where we regard
$\Spec(M)$ as a multiset where each eigenvalue of $M$ occurs with the
corresponding
multiplicity. It turns out, even in
the more general case where $(\Adj, \tau)$ is not a matrix algebra,
many of the general facts from linear algebra carry over.
\begin{lemma}[\cite{NicaS06}, Lecture 3]\label{lem:eig_gen}
    Let $(\Adj, \tau)$ be a $C^*$-probability space with a faithful
    trace $\tau$.
    \begin{enumerate}[(i)]
        \item For any $a \in \Adj$, $\Spec(a)$ is a compact subset of
            $\set*{z \in \CC : \Abs*{z} \le \norm*{a}}$.
        \item When $a \in \Adj$ is normal, i.e.~$a^* a = a a^*$, then
            $\norm*{a} = \max_{z \in \Spec(a)} |z|$ and $\supp \mu_a
            = \Spec(a)$.
        \item When $a \in \Adj$ is self-adjoint, i.e.~$a= a^*$, then
            $\Spec(a) \subset \RR$.
        \item When $a \in \Adj$ is positive, i.e.~$a = b^* b$ for
            some $b \in \Adj$, then $\Spec(a) \subset \RRp$.
    \end{enumerate}
\end{lemma}

\begin{defn}[Freeness]
    Let $(\Adj, \tau)$ be a $C^*$-probability space and
    $\Adj_1,\ldots, \Adj_r$ be sub-algebras of $\Adj$. We say that
    the algebras $\Adj_1,\ldots, \Adj_r$ are free from each other when we have
    $\tau(a_1\cdots a_n) = 0$ whenever
    \begin{enumerate}[(i)]
        \item $\tau(a_i) = 0$ for all $i \in [n]$,
        \item neighboring $a_i$ are from different sub-algebras,
            i.e., whenever $a_j \in \Adj_{i_j}$  for $j \in [r]$ then
            $i_j \ne i_{j + 1}$
    \end{enumerate}
\end{defn}
%
Finally, we recall the free-product of groups.
\begin{defn}[Free Product of Groups]
    Let $G_i$ for $i \in I$ be groups with identity elements $e_i \in
    I$ respectively. The free product $G = *_{i \in I} G_i$ is
    the group which are generated by all the elements of $G_i$ for $i
    \in I$, subject to
    \begin{enumerate}[(i)]
        \item the relations within each $G_i$ for $i \in I$,
        \item the identity $e_i \in G_i$ being identified with the
            identity $e \in G$.
    \end{enumerate}
\end{defn}
\begin{xmpl}[Examples of Free Products]
    For our purposes we will mostly deal with the free-group $\Fr_r$
    ($r \in \NN$) which are generated by $r$ generators $g_1,\ldots,
    g_r$ which satisfy no relations among themselves.  $\Fr_r$ can be
    identified with the free-product of $\ZZ$ with itself $r$-times,
    i.e., $\Fr_r \cong \ZZ^{* r}$.

    We will also deal with the group $\Gam_d = \ZZ_2^{* d}$, which is
    generated by $d$ generators $h_1,\ldots, h_d$ which are free from
    each other and are subject to the relations $h_i^2 = e$ for all
    $i \in [d]$, where $e \in \Gam_d$ is the identity element of $\Gam_d$.
\end{xmpl}

\subsection{Probability Transforms}
\label{sub:probability_transforms} 
Let $\mu$ be a probability measure $\mu$ on $\RR$, the Cauchy
transform $G_\mu$ of $\mu$ is defined by,
\begin{equation}
    G_\mu(z) = \int_\RR \frac{1}{z-t}~\dif\mu(t)\qquad\forall~z \in \CC^+
    \label{eq:cauchy_def}
\end{equation}
where $\CC^+ = \set*{ x + \imath \cdot y \col x, y \in \RR; y > 0 }$
is the complex upper half plane. Suppose $\mu$ is compactly supported
and write $\rad := \sup\set*{|t| \col t \in \supp(\mu)}$. Throughout,
we will think of $\mu$ as the spectral measure of a bounded element
in a non-commutative probability space.

Then, we have the following series expansion
\begin{equation}
    G_\mu(z) = \sum_{j = 0}^\infty
    \frac{m_j(\mu)}{z^{j+1}}\qquad\forall~z~\tst~|z| > \rad
    \label{eq:cauchy_series}
\end{equation}
where
$m_j(\mu) = \int_\RR t^j~\dif \mu(t)$ denotes the $j$-th moment of
the distribution $\mu$, and we set $m_0(\mu) = 1$. We note that when
$\mu := \mu_a$ is the spectral distribution of a self-adjoint element
$a \in \Adj$ in a $C^*$-probability space $(\Adj, \tau)$, i.e.~$a^* =
a$, we have
$m_j(\mu) = \tau(a^j)$.

We define the $M$-transform $M_\mu$ of the distribution
as,
\begin{equation}
    M_\mu(z) = \sum_{j = 1}^\infty m_j(\mu) \cdot z^j = z^{-1} \cdot
    G_\mu(z^{-1})-1\qquad\forall~z~\tst~|z| < \frac{1}{\rad}.
    \label{eq:emm_def}
\end{equation}
where the latter equality is obtained by manipulating the series in
\eqref{eq:cauchy_series}. We note that the radius of
convergence of the $M$-transform is precisely $\norm*{a}^{-1}$
when $\mu := \mu_a$ is the spectral measure of a self-adjoint element
$a \ne 0$ in a $C^*$-probability space $(\Adj, \tau)$ with a faithful
trace, by \cref{lem:norm_def} and the Cauchy-Hadamard Theorem
\cite[Theorem 2.4.1]{Wilf94}. That is to say, in this scenario the
parameter $\rad$, which
denotes the radius of convergence of the series expansion of the
Cauchy transform
\eqref{eq:cauchy_series}, occuring in
\eqref{eq:cauchy_series} and \eqref{eq:emm_def} corresponds to $\norm*{a}$.

We will write $\Minv_\mu$ for the inverse
of $M_\mu$ under
composition as a formal power series, i.e.,
$\Minv_\mu\parens*{M_\mu(z)} = z.$

Finally, we define the $S$-transform $S_\mu$ of the measure $\mu$
(satisfying $m_1(\mu) \ne 0$),
\begin{equation}
    S_\mu(z) = \frac{1+z}z \cdot \Minv_\mu(z)
    \label{eq:ess_transform}
\end{equation}
which we will regard as a formal power series. Since the support of
$\mu$ is bounded by $\rad$, it can be shown that $S_\mu$ is
analytic for all $z$ where $0 < |z| < \ee(\rad)$ for some small
enough $\ee(\rad)$.

We will now point out two properties of the $S$-transform. The first
one is rather straightforward,
\begin{fact}\label{fac:tracial_ess}
    Let $(\Aa, \tau)$ be a $C^*$-probability space, with
    normal $a, b \in \Aa$ -- $a^* a = a a^*$ and $b^* b = b b^*$ --
    satisfying $\tau(ab) \ne 0$. Defining the $S$-transforms of
    $S_{\mu_{ab}}$ and $S_{\mu_{ba}}$ in terms of the moment series,
    \eqref{eq:ess_transform}
    for the spectral laws $\mu_{ab}$ and $\mu_{ba}$ of $ab$ and $ba$
    respectively, we have
    \[ S_{\mu_{ab}}(z) = S_{\mu_{ba}}(z).\]
\end{fact}
\begin{proof}
    Notice that by the tracial property, we have
    \[ \tau( (ab)^p ) = \tau(a b (ab)^{p-1}) = \tau(b (ab)^{p-1} a) =
    \tau( (ba)^p )  \]
    for all $p \in \NN$. Thus, $m_p(\mu_{ab}) = m_p(\mu_{ba})$ for
    all $p \in \NN$, and consequently $\mu_{ab}$ and $\mu_{ba}$ have
    the same $M$-transform. The statement now follows from
    \eqref{eq:ess_transform}.
\end{proof}
The second property we will present is the main utility of the
$S$-transform, namely, that it reduces \emph{multiplicative free
convolution} to pointwise multiplication.
\begin{theorem}[Corollary~18.17, \cite{NicaS06}]\label{thm:mult_ess}
    Let $(\Aa, \tau)$ be a $C^*$-probability space, with
    $a, b \in \Aa$ being positive elements free from each other --
    i.e.~there exists $x, y \in \Aa$ such that $a = x^*x$ and $b =
    y^* y$ -- and satisfying $\tau(a), \tau(b) \ne 0$. Then, writing $\mu_c$ for
    the spectral measure of $c \in \Aa$, we have
    \[ S_{\mu_{ab}}(z) = S_{\mu_a}(z) \cdot S_{\mu_b}(z).\]
\end{theorem}
\begin{remark}
    Notice that even though $ab$ is not a normal element, we have
    \[ \tau((ab)^p) = \tau( (x^* x b)^p) = \tau(x^* (xb x^*)^{p-1} x
    b) = \tau( (x bx^*)^p).\]
    Thus, it has the same spectral distribution as a normal element.
\end{remark}
Finally, we introduce Pringsheim's theorem,
\begin{theorem}[Pringsheim's Theorem; Theorem 4.6,
    \cite{FlajoletS09}]\label{thm:pringsheim}
    If $F(z) = \sum_{n \ge 0} a_n z^n$ be a power series with
    non-negative real coefficients $a_n \in \RRp$ for all $n \ge 0$
    and radius of convergence $\rad$. Then, $F$ necessarily has a
    singularity at $z = \rad$.
\end{theorem}
\begin{remark}\label{rem:pringsheim}
    The relevance of Pringsheim theorem (\cref{thm:pringsheim}) is as
    follows: Let $(\Aa,\tau)$ be a $C^*$-probability space with a
    faithful trace, and $a\in\Aa$ be a positive non-zero element.

    Since the $M$-transform $M(z) := M_{\mu_a}(z)$,
    \cref{eq:emm_def}, of $\mu_a$
    has non-negative coefficients and
    $\tau(a^p)^{1/p}\to\norm*{a}$ as $p \to \infty$, its radius of convergence
    is $\norm*{a}^{-1}$. Pringsheim's theorem (\cref{thm:pringsheim})
    therefore implies
    that $z=\norm*{a}^{-1}$ is a singularity of $M$.

    This reduces the problem of analyzing the norm $\norm*{a}$ to a
    problem of understanding the singularities of $M(z).$
\end{remark}
\subsection{Kesten-McKay Measure}\label{sub:kesten_mckay_measure} 
For us, the measure of interest will be the (square) of the
Kesten-McKay measure: Let $g_1,\ldots, g_r \in \Fr_r = \ZZ^{* r}$ be the
generators of the free group (which is the free-product of $\ZZ$ with
itself $r$ times) and set $d = 2r$.
Then, the adjacency matrix of the infinite $d$-regular tree $\Tre_d$ can be
identified with,
\[ t = \lambda(g_1) + \lambda(g_1)^* + \cdots + \lambda(g_r) +
\lambda(g_r)^*.\]
Alternatively, writing $\Gam_d = \ZZ_2^{* d}$ for the free-product of
$\ZZ_2$ with itself $d$ times, and $h_1,\ldots, h_d$ free generators
for it we can identify
\[ t = \lambda(h_1) + \cdots + \lambda(h_d).\]
This formulation also allows odd $d \ge 2$.

By \cite{Kesten59}, it is known that
$\norm*{t} = 2\sqrt{d-1}$. Further, by \cite{McKay81} it is known that the
spectral measure of $t$ is given by the \emph{Kesten-McKay} measure,
\begin{equation}
    \dif\mu_t(x) =
    \begin{cases}
        \frac{d \sqrt{4 (d-1) - x^2}}{2\pi (d^2 - x^2)} &
        \tif~\Abs*{x} \le 2\sqrt{d-1},\\
        0 & \totw.
    \end{cases}
    \label{eq:km_def}
\end{equation}
The Cauchy transform of $\mu_t$ is given by
\[
    G_{\mu_t}(z) = \frac{(d-2)z - d\sqrt{z^2 - 4(d-1)}}{2(d^2 - z^2)}
\]
with the branch chosen such that $G_{\mu_t}(z) \sim z^{-1}$ at infinity,
cf.~\cite[Claim~3.3]{CohenM23}. We can compute $G_{\mu_{t^2}}(z) =
G_{\mu_t}(\sqrt{z})/\sqrt{z}$ as in \cite[Claim~3.5]{CohenM23}, and
arrive at \cite[Claim~3.4(1)]{CohenM23},
\[
    G_{\mu_{t^2}}(z) = \frac{d\sqrt{z - 4(d-1)} -
    (d-2)\sqrt{z}}{2\sqrt{z}(z - d^2)}
\]
Similarly, using \eqref{eq:emm_def}, i.e., $M_{\mu_{t^2}}(z) = z^{-1}
G_{\mu_{t^2}}(z^{-1})- 1$, it follows that
\[
    M_{\mu_{t^2}}(z) = \frac{d\sqrt{1 - 4(d-1)z} - (d-2)}{2(1 - d^2 z)}-1.
\]
We can similarly compute as in \cite[Claim 3.4(2)]{CohenM23}
\begin{equation}
    \Minv_{\mu_{t^2}}(z) = \frac{z(d+z)}{d^2 (1+z)^2}.
    \label{eq:t2_minv}
\end{equation}
Combining \eqref{eq:t2_minv} with the definition
\eqref{eq:ess_transform} of the $S$-transform, we obtain
\begin{equation}
    S_{\mu_{t^2}}(z) = \frac{d+z}{d^2 (1+z)}.
    \label{eq:t2_ess}
\end{equation}
\subsection{Weak and Strong
Convergence}\label{sub:weak_and_strong_convergence} 
Let $\rv X^N = (X_1^N, \ldots, X_r^N)$ be a sequence of random
tuples of matrices and let $\rv x=  (x_1,\ldots, x_r)$ be
a sequence of elements in a $C^*$-probability space $(\Adj,
\tau)$. We say that $\rv X^N$ weakly converges to $\rv x$, if
\begin{equation}
    \Pp\sqbr*{\Abs*{\tr\parens*{P(\rv X^N)} - \tau(P(\rv
    x))} \ge \ee} = o_N(1)
    \label{eq:weak_convg_def}
\end{equation}
for every non-commutative $*$-polynomial $P \in \CC^*\langle
x_1,\ldots, x_r\rangle$ and constant $\ee > 0$.

Similarly, we say $\rv X^N$ strongly converges to $\rv x$, if
\begin{equation} \Pp\sqbr*{\Abs*{\norm*{P(\rv X^N)} - \norm*{P(\rv
    x)}} \ge \ee} = o_N(1)
    \label{eq:strong_convg_def}
\end{equation}
for every non-commutative $*$-polynomial $P \in \CC^*\langle
x_1,\ldots, x_r\rangle$ and constant $\ee > 0$.

We recall the following consequence of weak convergence,
\cite[Lemma~2.13(c)]{vanHandel26}
\begin{lemma}\label{lem:weak_bd}
    Suppose $\rv X^N = (X_1^N, \ldots, X_r^N)$ weakly converges to
    $\rv x = (x_1,\ldots, x_r)$ in some $C^*$-probability
    space $(\Adj, \tau)$ with a faithful trace $\tau$. Then,
    \[ \norm*{P(\rv X^N)} \ge \norm*{P(\rv x)} -
    o(1)\quad\textrm{with probability}\quad 1 - o(1).\]
\end{lemma}


\subsection{Strong Convergence Phenomenon for Permutation Matrices}
\label{sub:strong_convergence_phenomenon} 
We recall the results of strong convergence for permutation matrices
proven by \cite{BordenaveC19} and \cite{ChenGVTvH26}. Throughout, let
$g_1,\ldots, g_r$ be the canonical generators freely generating
$\Fr_r$. We will set,
\begin{equation}\label{eq:u_def}
    u_i = \lambda(g_i)\qquad\textrm{for
    all}\qquad i = 1,\ldots, r
\end{equation}
where $\lambda(\bullet)$ denotes the left-regular representation of
$\Fr_r$, i.e., $\lambda(g) \delta_w = \delta_{g w}$ for all $g, w \in
\Fr_r$.

Similarly, we will let $h_1,\ldots, h_s \in \Gam_s = \ZZ_2^{* s}$ be
the canonical involutive elements freely generating $\Gam_s$
\begin{equation}\label{eq:m_def}
    m_i = \lambda(h_i)\qquad\textrm{for
    all}\qquad i = 1,\ldots, s.
\end{equation}
We introduce the shorthand $\rv u = (u_1,\ldots, u_r)$ and $\rv m =
(m_1, \ldots m_s)$ and assume that $\rv u$ is free from $\rv m$,
i.e., $\rv u$ and $\rv m$ together freely generate $\Fr_r * \
\Gam_s$. For a
sequence of matrix $A^N$ we will write $\left. A^N\right|_{\perp 1}$
for its projection away from the subsace $\sspan(\One)$, i.e.~$\left.
A^N\right|_{\perp 1} = (I - J^N)A^N (I - J^N)$ where $J^N = \One\One^\top/N$.

We first recall the following asymptotic strong convergence result by
\citeauthor{BordenaveC19},
\begin{theorem}[Asymptotic Convergence,
    \cite{BordenaveC19}]\label{thm:asymptotic_bc}
    Let $D, r, s\in \NN$ be constants. For any $N \in \NN$, we let
    $\rv U^N = (U_1^N, \ldots, U_r^N)$ be a
    sequence of iid uniformly random
    permutation matrices, and $P \in M_D(\CC) \tensor \CC^*\langle
    x_1,\ldots, x_r \rangle$ a non-commutative $*$-polynomial with matrix
    coefficients. Then, for all $\ee > 0$, we have
    \[\Pp\sqbr*{
            \Abs*{\norm*{\left. P\parens*{ \rv U^N }\right|_{\perp
                1}} - \norm*{P(\rv
        u)}_{\Cred(\Fr_{r})} } \ge \ee } = o_N(1)
    \]
    where $\rv u$ is as defined in \eqref{eq:u_def}.

    When $N$ is even, writing $\rv M^N = (M_1^N, \ldots, M_s^N)$ for a tuple of
    iid uniformly random matching matrices\footnote{these are
        permutation matrices belonging to involutions $\sigma \in
    \Sym_N$, i.e.~$(M_j^N)^2 = I^N$ for all $j = 1,\ldots s$}
    we have that for any
    $P \in M_D(\CC) \tensor
    \CC^*\langle x_1,\ldots, x_r, y_1, \ldots, y_s\rangle$
    \begin{equation}
        \Pp\sqbr*{
            \Abs*{\norm*{\left. P\parens*{ \rv U^N ,
                \rv M^N }\right|_{\perp 1}} - \norm*{P(\rv
        u, \rv m)}_{\Cred(\Fr_{r} * \Gam_s)} } \ge \ee } = o_N(1)
    \end{equation}
    where $\rv m$ is as defined as \eqref{eq:m_def} and $\rv u$ is as
    defined as \eqref{eq:u_def}.
\end{theorem}
In other words, \cref{thm:asymptotic_bc} states that as $N$ tends to
infinity, the operator norm of any $*$-polynomial of fixed degree
with matrix coefficients in iid uniform permutation or matching
matrices is approximated by the norm of the same
polynomial applied to the tuple $\rv u$ which generates the reduced group
$C^*$-algebra of $\Fr_r$ or $\Fr_r * \Gam_s$.

Recently \citeauthor{ChenGVTvH26} has proven a non-asymptotic version
of \cref{thm:asymptotic_bc},
\begin{theorem}[Non-Asymptotic Convergence,
    \cite{ChenGVTvH26}]\label{thm:nonasymptotic}
    Let $D \in \NN$, $\rv U^N = (U_1^N, \ldots, U_r^N)$ be a
    sequence of iid uniformly random
    permutation matrices ($r \ge 2$), and $P \in M_D(\CC) \tensor
    \CC^*\langle
    x_1,\ldots, x_r \rangle$ a non-commutative and self-adjoint
    $*$-polynomial with matrix
    coefficients of degree $q_0$. Then,
    \[ \Pp\sqbr*{ \norm*{\left. P\parens*{\rv
            U^N}\right|_{\perp 1}}  \ge \norm*{P(\rv
        u)}_{\Cred(\Fr_r)}+ \ee} \lesssim \frac{D}{N} \cdot
        \parens*{ \frac{K q_0
    \log r}{\ee} }^8 \cdot \log\parens*{\frac{ e K}{\ee} } \]
    where $K = \norm*{P}_{M_D(\CC) \tensor C^*(\Fr_r) }$ and $0 \le
    \ee < K - \norm*{P(\rv u)}$.
\end{theorem}
\begin{remark}\label{rem:c*norm}
    Suppose $D = 1$, we note that the norms
    $\norm*{\bullet}_{\Cred(\Fr_r)}$ and
    $\norm*{\bullet}_{C^*(\Fr_r)}$ are not the same. Nevertheless, we
    can still bound the later norm by taking a sum over the absolute
    value of coefficients of $P$.
\end{remark}

We will also make use of the following constructivization of
\cref{thm:asymptotic_bc} by
\citeauthor{OdonnellW20},
\begin{theorem}[\cite{OdonnellW20}, Simplified Version of Theorem
    10.13]\label{thm:constructive_convg}
    Fix constants $R, \ee > 0$ and let $N \in \NN$. Let $\mathcal{P}
    \subset \CC^*\langle x_1,\ldots, x_r,
    y_1,\ldots, y_s\rangle$ be a collection of $*$-polynomials, such
    that every polynomial in $\mathcal P$
    has degree at most $k$ and the coefficients of polynomials in
    $\mathcal P$ are bounded by $R$.

    On input $N$, there is an algorithm which runs in time $N^{O_{k,
    R, r, s, \ee}(1)}$ and outputs a sequence of
    (deterministic) permutation matrices $\rv U^{N'} = (U_1^{N'}, \ldots,
    U_r^{N'})$ and (deterministic) matching matrices $\rv M^{N'}
    = (M_1^{N'},
    \ldots, M_s^{N'})$ such that simultaneously for every $P \in
    \mathcal P$, we have
    \[ \Abs*{\norm*{\left. P\parens*{ \rv U^{N'},  \rv
            M^{N'}}\right|_{\perp 1} }  -\norm*{P(\rv u, \rv
        m)}_{C^*_{\mathrm{red}}(\Fr_r * \Gam_s) }} \le \ee\qquad\forall~P
    \in \mathcal P\]
    where $N \le N' \le (1 + o(1)) N$.
\end{theorem}

\section{Rotating Infinite Trees and Expanders}
\label{sec:rotating_infinite_trees_and_expanders} 
\subsection{Overview}\label{sub:overview_rot} 
Let $\rv V^N = (V_1^N, \ldots, V_{ k-1 }^N)$ be a sequence of
iid permutation matrices. \citeauthor{CohenM23} \cite{CohenM23}
define the rotating random walk over an $(N, d, \lambda)$-graph $G$
with adjacency matrix $A^N$ as,
\begin{equation}\label{eq:rot_walk_def}
    \Arot(\rv V^N) = A^N V_1^N A^N V_2^N \cdots V_{k-1}^N A^N \cdot
    A^N V_{k-1}^{N*} A^N \cdots V_2^{N*} A^N V_1^{N*} A^N.
\end{equation}
For simplicity, we will consider their result in the case where $G$
is Ramanujan and has high girth, i.e, $\lambda \le 2\sqrt{d-1}$ and
$\girth(G) = \Theta(\log_{d-1}(n))$. They show that there
exists a \emph{deterministic} sequence of permutation matrices $\rv \Pi^N
= (\Pi_1^N, \ldots, \Pi_{ k-1 }^N)$ such that
\[ \lambda( \Arot(\rv \Pi^N) ) \le \parens*{1 + \frac{1}{k}}^k (k+1)
d^k + o(1)\]
so long as $k \log (kd) = o(\log g)$ with a similar extension for
non-Ramanujan graphs of high girth after an initial burn in period,
\cite[Theorems~7.3-7.4]{CohenM23}. In probabilistic
language, their result can be thought as showing
\begin{equation}
    \Pp_{\rv V^N}\sqbr*{ \lambda(\Arot(\rv V^N)) \ge \parens*{1 +
        \frac{1}{k}}^k (k+1)
    d^k + \ee} < 1
    \label{eq:cm_existential}
\end{equation}
for all \emph{fixed} $\ee > 0$ where $G := G_N$ is taken from a
family of high-girth Ramanujan graphs of increasing sizes, as $N \to
\infty$, say as in \cite{LubotzkyPS86, Morgenstern94}. In
\cite{CohenCMP25}, \citeauthor{CohenCMP25} show that this bound is
\emph{essentially} tight. Namely, for \emph{any} sequence of
deterministic permutation
matrices $\rv \Pi^N = (\Pi^N_1, \ldots, \Pi_{k-1}^N)$, we have
\begin{equation}
    \lambda\parens*{\Arot\parens*{\rv \Pi^N}} \ge
    \frac{(k+1)^k}{k^k} \cdot (k+1 )\cdot (d-1) d^{k-1} - o_N(1).
    \label{eq:lower bound_fc}
\end{equation}
So up to the lower order $o_N(1)$ error term, the result of
\cite{CohenM23} is tight.

The main objective of this section is to show that when $G$ is a
random graph or a suitable base graph, the RHS of
\eqref{eq:cm_existential} can be taken to be $o_N(1)$, using the
non-asymptotic convergence result
\cref{thm:nonasymptotic} due \citeauthor{ChenGVTvH26}. We note that
this confirms an observation made in \cite{CohenM23}, that their
bound seems to capture the \emph{typical} case behaviour of \emph{rotating
random walks}.

Let $g_1,\ldots, g_r, h_1,\ldots, h_{k-1}$ the generators of $\Fr_{r
+ k -1}$. We will write,
\[ u_j = \lambda(g_j),~v_i = \lambda(h_i)\qquad\forall~j \in [r], i
\in [k-1]\]
where $\lambda(\bullet)$ denotes the left-regular representation of
$\Fr_{r + k - 1}$. As in \cref{sub:kesten_mckay_measure} we introduce
the shorthand
\[ t =  u_1 + u_1^* + \cdots + u_r + u_r^*\]
for the adjacency matrix of the $d = (2r)$-regular infinite tree
$\Tre_d$.  We observe that the limiting object for the (random)
rotating walk
on a random graph, is the \emph{rotating infinite tree} $a$
defined by setting,
\begin{equation}
    a = t v_1 t v_2 t \cdots v_{k-1} t^2 v_{k-1}^* \cdots v_2^* t v_1^* t
    \label{eq:a_def}
\end{equation}
Notice that
$a = P(\rv u, \rv v)$ for a non-commutative $*$-polynomial $P \in
\CC^*\langle x_1,\ldots, x_r, y_1, \ldots, y_{k-1} \rangle$ of degree
$4k$. Further, applying the same polynomial to $\rv U^N =
(U_1^N,\ldots, U^N_r)$ and $\rv V^N = (V_1^N, \ldots, V_{k-1}^N)$ yields,
\[
    P(\rv U^N, \rv V^N) = \Arotr\parens*{\rv U^N, \rv
    V^N}\qquad\textrm{where}\qquad \Arand(\rv U^N) = \sum_{j = 1}^r
    \parens*{U_j^N +
U_j^{N*}}.\]
In words, $P(\rv U^N, \rv V^N)$ describes the rotating random walk
with random permutations $\rv V^N$ of a random $d = (2r)$-regular
graph with adjacency matrix $\Arand := \Arand(\rv U^N)$ drawn from the
permutation model.
In particular \cref{thm:asymptotic_bc} implies,
\begin{equation}
    \Pp_{\rv U^N, \rv V^N}\sqbr*{ \lambda\parens*{\Arotr\parens*{\rv U^N,
        \rv V^n }} \ge
    \norm*{a}_{\Cred(\Fr_{r+k-1})} + \ee} = o_N(1) \qquad\forall~\ee > 0
    \label{eq:showcase_tree}
\end{equation}
in fact as we will see the $o_N(1)$ term can be bounded by above
using the non-asymptotic bounds of \cite{ChenGVTvH26}
(\cref{thm:nonasymptotic}). We note that when $N$ is even
\eqref{eq:showcase_tree}
remains valid when the $d$-regular graph
$\Arand(\rv M^N)$ from the matching model, i.e.
$\Arand(\rv M^N) = \sum_{i = 1}^d M_i^N$. Since we cannot make this
result non-asymptotic, we will not prove this explicitly.

\paragraph{Proof Idea.} We will analyze the norm
$\norm*{a}_{\Cred(\Fr_{r + k - 1)}}$ of $a$
from \eqref{eq:a_def} by
using Pringsheim's theorem (\cref{thm:pringsheim,}),
i.e., by locating the first singularity of the $M$-transform $M(z):=
M_{\mu_a}(z)$ of $\mu_a$ from \eqref{eq:emm_def} on the positive real
axis which we know corresponds to $\norm*{a}^{-1}_{\Cred(\Fr_{r+ k-1})}$
by \cref{rem:pringsheim}.

To get a handle on
the $M$-transform, we will argue that $a$ can be thought as a product
of trees which are free from each other, and use the $S$-transform
\eqref{eq:ess_transform} of the (squared) Kesten-McKay measure
\eqref{eq:t2_ess} and the multiplicative property of $S$-transforms
\cref{thm:mult_ess} to obtain the $S$-transform of $a$. Since the
$S$-transform and the inverse
$M$-transform are related,
this will allow us to argue about $M$-transform by means of the
inverse $M$-transform $\Minv(z)$. It turns out by using analytic
ideas, the task of locating the first singularity of $M$ can be
reduced to finding a critical point of $\Minv(z)$. This will ensure
that the only
distributional input we require to compute the norm of $a$ \emph{exactly} is
the $S$-transform
of the (squared) Kesten-McKay measure of the $(2r)$-regular tree,
which we readily know \eqref{eq:t2_ess}.
\subsection{Bounding the Norm of an Infinite Rotating Tree}
\label{sub:bounding_the_norm_of_an_infinite_rotating_tree} 
Our main observation in this section is that although the definition
of the rotating infinite tree $a$ in \eqref{eq:a_def} has the same
tree operator $t$ used over
and over, the interleaving with $v_i$-s allows us to
\enquote{decouple} successive occurences of $t$. Below, we will
formalize this idea.

We define the words,
\begin{equation}
    w_j = h_1 h_2 \cdots h_{j}\qquad\forall~j = 1,\ldots, k-1
    \label{eq:doubleu}
\end{equation}
and set $w_0 = e$ (empty word in $\Fr_{r + k - 1}$).
We also define the conjugated trees,
\begin{equation}
    t_j = \lambda(w_j) \cdot t \cdot \lambda(w_j)^*\qquad\forall~j =
    0, \ldots, k-1.\label{eq:tj_def}
\end{equation}

\begin{obsv}\label{obsv:freeness_tj}
    Let
    \[ \Gr_j = w_j \cdot \langle g_1,\ldots, g_r \rangle \cdot w_j^{-1}
    \le \Fr_{r + k -1}\qquad\forall~j = 0, \ldots, k- 1.\]
    Then, the subgroups $\Gr_0,\ldots, \Gr_{k-1}$ of
    $~\Fr_{r + k}$ are free from each other. Consequently,
    $t_0,\ldots, t_{k-1}$ from \eqref{eq:tj_def} are free from each other
    with the same spectral
    distribution $\mu_t$ from \eqref{eq:km_def}.
\end{obsv}
\begin{proof}
    Let $i_1,\ldots, i_\ell \in [k]$ be given such that $i_1 \ne i_2$,
    $i_2\ne i_3$, \ldots, $i_{\ell - 1} \ne i_\ell$. Similarly, let
    $\gamma_1,\ldots,
    \gamma_\ell$ be words such that $\gamma_j \in \Gr_{i_j} \setminus
    \set*{e}$ for each
    $j \in [\ell]$. Notice that we can write,
    \[ \gamma_1 \cdots \gamma_\ell = w_{i_1} \bar\gamma_{i_1}
        w_{i_1}^{-1} \cdots w_{i_j}^{-1} w_{i_{j+1}} \bar\gamma_{i_{j+1}}
        w^{-1}_{i_{j+1}} w_{i_{j+2}} \cdots w_{i_{\ell-1}}^{-1} w_{i_\ell}
    \bar\gamma_{i_\ell} w_{i_\ell}\]
    where each $\bar\gamma_{i_j} \in \langle g_1,\ldots, g_r
    \rangle$. Notice that since the $w_j$-s from \eqref{eq:doubleu} are
    free from $\bar
    \gamma$-s, the only possible places for reduction are the
    connection points, and it is evident that we cannot have
    $w_{j}^{-1} w_k =e$ unless $j = k$ . Consequently, we have
    $\gamma_1 \cdots \gamma_\ell \ne e$. Establishing the
    freeness of subgroups.

    Recall that the canonical trace
    of $\Cred(\Gr)$ is defined by,
    \[ \tau\parens*{ \sum_{g \in \Gr} \alpha_g \cdot
    \lambda(g)} = \alpha_e\]
    for any sequence of $\alpha_g$-s such that only finitely many
    $\alpha_g$-s are non-zero. Thus, the freeness of $\Gr_1,\ldots,
    \Gr_k$ immediately implies the freeness of $t_i$. To see that the
    spectral measure of these coincide with the spectral measure of
    $\Tre_{d}$, we note
    \[ \tau\parens*{ t^{p} } = \tau\parens*{ t_j^p }\qquad\forall~j =
    1,\ldots, k; p \in \NN.\]
    Thus, the theorem follows by the definition of the spectral measure
    (\cref{defn:spec_def}).
\end{proof}
The utility of $t_j$-s from \eqref{eq:tj_def} and
\cref{obsv:freeness_tj} is due to the
following identity concerning the infinite rotating tree from
\eqref{eq:a_def},
\begin{equation}
    a = (t_0 \cdots t_{k-1}) (t_0 \cdots t_{k-1})^*
    \label{eq:a_alternative_def}
\end{equation}
which is obtained by observing $w_{j-1}^{-1} w_j = h_j$ for each $j =
1,\ldots, k-1$.
\begin{lemma}\label{lem:main_lemma}
    Let $a$ be as in \eqref{eq:a_def}, writing $\mu_a$ and
    $\mu_{t^2}$ for the spectral laws of $a$ and $t^2$, we have
    \[ S_{\mu_a}(z) = \parens*{S_{\mu_{t^2}}(z)}^k =
    \parens*{\frac{d+z}{d^2 (1+z)} }^k.\]
\end{lemma}
\begin{proof}
    For simplicity, we will write $S_c$ to denote the $S$-transform
    of the spectral measure $\mu_c$ in place of $S_{\mu_c}$.

    We have,
    \begin{align*}
        S_{a}(z) &~=~S_{(t_0 \cdots t_{k-1})(t_0 \cdots t_{k-1})^*}(z),\\
        \textrm{(by associativity)}&~=~
        S_{t_0 (t_1 \cdots t_{k-1}) (t_1 \cdots t_{k-1})^* t_0 }(z),\\
        \textrm{(by \cref{fac:tracial_ess})}&=~S_{t_0^2 (t_1 \cdots
        t_{k-1})(t_1 \cdots t_{k-1})^*}(z),\\
        \textrm{(by \cref{thm:mult_ess})}&~=~S_{t_0^2}(z) \cdot
        S_{(t_1 \cdots t_{k-1}) (t_1 \cdots t_{k-1})^*}(z)
    \end{align*}
    where the freeness assumption needed for \cref{thm:mult_ess} is
    asserted by \cref{obsv:freeness_tj}. Now, repeating the same argument
    $(k-1)$ times more we obtain
    \[ S_a(z) = S_{t_0^2}(z) \cdots S_{t_{k-1}^2}(z) =
    \parens*{\frac{d+z}{d^2 (1+z)} }^k\]
    since all $t^2_i$-s have the same spectral measure. The result now
    follows from \eqref{eq:t2_ess}.
\end{proof}
\begin{theorem}[Norm of the Infinite Rotating Tree]\label{thm:inf_rot_tree}
    Let $a$ be as in \eqref{eq:a_def}. Then,
    \[ \norm{a}_{\Cred(\Fr_{k+r-1})} = \frac{(k+1)^{k+1}}{k^k} \cdot
    (d-1) d^{k-1}. \]
\end{theorem}
\begin{proof}
    By \cref{lem:main_lemma} and \eqref{eq:ess_transform}, we have
    that the \emph{inverse} $M$-transform of $\mu_a$ is given by
    \begin{equation}
        \Minv(z) = \frac{z}{1+z} \cdot S_a(z) = \frac{z}{1+z} \cdot
        \parens*{\frac{d+z}{d^2 (1+z)}}^k = \frac{z (d+z)^k}{d^{2k}
        (1+z)^{k+1}}.
        \label{eq:complicated_minv}
    \end{equation}
    The advantage of working with the inverse $M$-transform is evident from
    \cref{eq:complicated_minv}: Even though $M$ itself is not known
    to us, we can get a simple rational expression for
    its inverse. Starting at the origin, we can recover $M$ by
    analytically inverting the increasing branch of $I(z) = \Minv(z)$, as long
    as $I'$ does not vanish. We will show that this branch ends
    at a first positive critical point $z^*$, where analytic
    inversion becomes impossible. Since Pringsheim's theorem rules
    out a smaller radius of convergence for $M$
    (\cref{thm:pringsheim,rem:pringsheim}),
    the critical value $I(z^*)$ is precisely its radius of
    convergence, and $\norm*{a}=I(z^*)^{-1}$.

    We now carry out this plan: For $y > \norm*{a}_{\Cred(\Fr_{r+k-1})} =
    \sup~\supp\parens*{\mu_a}$ (where the final equality is due to
    \cref{lem:eig_gen}), \eqref{eq:complicated_minv} and the
    relation $\Minv(M(y^{-1})) = y^{-1}$ implies
    that we can write
    \[ y= \frac{1}{\Minv(z)} = \frac{(1+z)^{k+1} d^{2k}}{z
    (d+z)^k}\qquad\textrm{where}\qquad z = M(y^{-1})\]
    where we note that $M(y^{-1})$ is well-defined as $y^{-1} <
    1/\norm*{a}_{{\Cred(\Fr_{r + k - 1})}}$.

    Writing $I(z) = \Minv(z)$ for convenience, we note that
    \[ \frac{I'(z)}{I(z)}= \frac1z + \frac k{d+z} - \frac{k+1}{1+z}.\]
    Thus, for $z^* = \frac{d}{k(d-1) - 1}$ we have $I'(z^*) = 0$. We
    can verify that this is a critical point since, $I'(z) > 0$ for
    $z \in [0, z^*)$ and $I'(z) < 0$ for $z > z^*$. In particular, we
    see that $I(z) = \Minv(z)$ is strictly increasing in the interval
    $[0, z^*)$ and its inverse $M(z)$ is analytic in $[0, \Minv(z^*))$.

    Let $\rad^{-1}$ be the radius of convergence of $M(z)$. By
    \cref{lem:eig_gen} and \eqref{eq:emm_def}, this means that
    $\rad = \norm*{a} = \sup \supp \mu_a$. By Pringsheim's
    theorem (q.v.~\cref{thm:pringsheim,rem:pringsheim}), since $M(z)$ is
    analytic in $[0,
    \Minv(z^*))$, we necessarily have $\rad^{-1} \ge \Minv(z^*)$.
    Assume by contradiction that
    $\rad^{-1} > \Minv(z^*) = I(z^*)$, then $M$ is analytic at
    $\Minv(z^*) = I(z^*)$
    and thus $M(\Minv(z^*)) =
    z^*$. Differentiating both sides, this would mean by the chain rule
    \[ 0 = I'(z^*) \cdot M'(I(z^*)) = 1 \]
    since we have picked $z^*$ to make the first derivative of
    $\Minv(z)$ vanish. Thus, we necessarily have $\rad^{-1} = \Minv(z^*)$.

    A direct computation using the formula
    \eqref{eq:complicated_minv} yields the result.
\end{proof}
\subsection{Rotating Walks on Random Graphs}
\label{sub:rotating_walks_on_random_graphs} 
We note that \cref{thm:inf_rot_tree,thm:asymptotic_bc}
readily imply that for fixed $k$ and $r$, we have that the rotating
random walk $\Arotr(\rv U^N, \rv V^N)$ over the random graph
$\Arand(\rv U^N)$ from the permutation model with the rotations $\rv
V^N$ \emph{almost} meets
the lower bound in \cite{CohenCMP25} except with vanishing
probability. Formally,
{\small
    \begin{equation}
        \label{eq:asymptotic_randrot}
        \Pp_{\rv U^N, \rv V^N}\sqbr*{
            \lambda\parens*{\Arotr\parens*{ \rv U^N, \rv
        V^N }} \ge  \frac{(k+1)^{k+1}}{k^k} \cdot (d-1) d^{k-1 } +  \ee } =
        o_N(1)\quad\textrm{where}\quad \Arand(\rv U^N) = \sum_{i
        =1}^r (U_i^N + U_i^{ N* }).
    \end{equation}
}
Now, we will appeal to \cref{thm:nonasymptotic} to show that we can
let $k$ and $d$ grow with the number of vertices $N$.
\begin{theorem}\label{thm:permutation_model_nonasymp}
    For every $\ee$ such that
    \[ 0 < \ee < d^{2k} - \frac{(k+1)^{k+1}}{k^k} \cdot (d-1) d^{k-1}\]
    we have
    \[
        \Pp_{\rv U^N, \rv V^N}\sqbr*{
            {\lambda\parens*{\Arotr\parens*{\rv U^N, \rv
        V^N }}\ge \frac{(k+1)^{k+1}}{k^k} \cdot (d-1) d^{k-1 }} + \ee }
        \lesssim \frac{1}{N} \cdot \parens*{\frac{d^{2k} 4k \cdot
        \log(r+k)}{\ee} }^8 \cdot \log\frac{e d^{2k}}{\ee}
    \]
    where $d = 2r$ and $\Arand := \Arand(\rv U^N)$ is a random $d$-regular
    graph drawn from the permutation model.
\end{theorem}
\begin{proof}
    Let $P$ be as in \cref{sub:overview_rot}, i.e., $P(\rv u, \rv v)
    = a$ and $a$ is as defined in \eqref{eq:a_def}. Notice that $P$ is
    a polynomial of degree $4k$ with $d^{2k}
    = (2r)^{ 2k }$ monomials each having the coefficient 1, and consequently
    $\norm*{P}_{C^*(\Fr_r)} \le d^{ 2k }$. Thus the statement
    immediately follows by \cref{thm:nonasymptotic}, since we have
    $\lambda(\Arotr(\rv U^N, \rv V^N)) = \norm*{\left. P(\rv U^N, \rv
    V^N)\right|_{\perp 1}}.$
\end{proof}
\subsection{Finding a Rotating Almost Ramanujan
Graph}\label{sub:finding_a_rotating_almost_ramanujan_graph} 
We observe that \cref{thm:inf_rot_tree} can be used in conjunction
with \cref{thm:constructive_convg} to find explicit concrete families
of graphs $G_N = (V_N, E_N)$ which are nearly Ramanujan, i.e.,
$\lambda(G_N) \le 2\sqrt{d-1} + \ee$ (for some bounded $\ee > 0$)
along side a family of sequences
of permutation matrices $\rv \Pi^N = (\Pi_1^N, \ldots, \Pi_{k-1}^N)$
such that $\Arot(\rv \Pi^N)$ almost saturates the bound in
\eqref{eq:lower bound_fc} up to an error of $\ee$.

\begin{corollary}\label{cor:rot_explicit}
    For each $N, r, k \in \NN$, and $\ee > 0$ there exists some $N'
    \in [N, (1+ o_N(1))
    N]$ such that there exists a $d = (2r)$-regular graph $G_{N'} =
    (V_{N'}, E_{N'})$ with adjacency matrix $A^{N'}$ and a sequence
    of permutation matrices $\rv
    \Pi^{N'} = (\Pi^{N'}_0,\ldots, \Pi^{N'}_{k-1})$ such that
    \begin{enumerate}[(i)]
        \item $\lambda(A_{N'}) \le 2\sqrt{d- 1} + \ee$,
        \item $\lambda\parens*{\Arot(\rv \Pi^{N'})} \le
            \parens*{\frac{k+1}k}^k (k+1) (d-1)d^{k-1} + \ee.$
    \end{enumerate}
    We recall that,
    \[\Arot(\rv \Pi^{{N'}'}) = A^{N'} \Pi_1^{N'} A^{N'} \Pi_2^{N'}
        \cdots \Pi_{k-1}^{N'} A^{N'} \cdot
        A^{N'} \Pi_{k-1}^{N'*} A^{N'} \cdots \Pi_2^{N'*}
        A^{N'} \Pi_1^{N'*} A^{N'}.
    \]
    Moreover, there is a deterministic algorithm which can compute
    the graph $G_{N'}$ and the permutations $\rv \Pi^{N'}$ in time
    $N^{O_{r, k, \ee}(1)}$.
\end{corollary}
\begin{proof}
    This is an immediate consequence of \cref{eq:a_def} and the
    standard result $\norm*{t}_{\Cred(\Fr_r)} \le 2\sqrt{d-1}$,
    qv.~\cref{sub:kesten_mckay_measure} --- where $t$ is the
    adjacency operator of the infinite $d := (2r)$-regular tree $\Tre_d$.

    Since both $a$ (as defined in \eqref{eq:a_def}) and $t = \sum_{i
    = 1}^{r} (u_i + u_i^*)$ can be described as $*$-polynomials of
    degree $O_{r, k}(1)$ in
    $\rv u = (u_1,\ldots, u_r)$ and $\rv v = (v_1, \ldots, v_{k-1})$
    where $u_i = \lambda(g_i), v_j =\lambda(h_j)$ for $i \in [r]$ and
    $j \in [k-1]$ where $g_1,\ldots, g_r, h_1,\ldots, h_{k-1}$ are
    the canonical free
    generators of $\Fr_{r + k-  1}$. Since the largest coefficient of
    these polynomials is $1$, the statement follows from
    \cref{thm:constructive_convg}.
\end{proof}


\section{Average Behavior of Derandomized Squaring}
\label{sec:average_behavior_of_derandomized_squaring} 
We equip $G$ with an
out-labelling $\Out_G(v, a)$ which is a bijective mapping between the
vertex neighborhood $N_G(v)$ of $v$ and $[d]$. We say this mapping is
consistent if every \emph{incoming edge} to a vertex $v \in V$ has
different labels. Notice that when $G = (V, E)$ is a $d$-regular
graph, by Birkhoff's theorem \cite{Birkhoff46}, we can decompose the
adjacency matrix
of $G$ into a sum of
$d$-permutation matrices,
\begin{equation}
    A_G = \sum_{j= 1}^d \Pi_j\qquad\textrm{for some
    choice of}\qquad
    \pi_1,\ldots, \pi_d \in \Sym_{V}
    \label{eq:a_decomp}
\end{equation}
where $\Pi_j$ denotes the natural permutation matrix underlying $\pi_j$
and we define $\Out_G(v, j) = \pi_j(v)$. This will be
the \emph{canonical} labelling we will assume throughout, and is
consistent by construction.

We recall now the derandomized square defined by
\citeauthor{RozenmanV05} \cite{RozenmanV05}: The
derandomized square $S = G \dsq H$ of the $(N, d)$-graph $G$ via the
$k$-regular graph $H = ([d], E')$ is an $(N, dk)$-graph on
the vertex set $V$, whose out-labelling is defined by
\[ \Out_S(v, (a, b)) = \Out_G( \Out_G(v, a), \Out_H(a,
b))\qquad\forall~a \in [d], b \in [k].\]
Assuming the labellings of $G$ and $H$ are consistent, this yields a
consistent labelling \cite[Proposition~4.2]{RozenmanV05}. The nomenclature
of \emph{derandomized} square is justified by the fact that $G \dsq
H$ can be thought as a spectral approximation of $G^2$. Indeed,
\citeauthor{RozenmanV05} show,
\begin{equation}
    \omega\parens*{G \dsq H} \le 1 - (1- \omega(G)^2)(1 - \omega(H))
    \le \omega(G)^2 + \omega(H)
    \label{eq:dsq_expansion}
\end{equation}
cf.~\cite[Theorem~4.4]{RozenmanV05}. In words, \eqref{eq:dsq_expansion}
states that when $H$ is a sufficiently good expander, $G \dsq H$ has
expansion comparable to $\omega(G^2) = \omega(G)^2$ despite the
former being of degree $d \cdot k$ as opposed to $d^2$. In
\cite{CohenCMP25}, \citeauthor{CohenCMP25} prove a lower bound for
the expansion of $G \dsq H$ when $H$ is vertex-transitive.
Concretely, they show
\begin{equation}
    \lambda(G \dsq H) \ge \Lambda_H - o(1)\qquad\textrm{where}\qquad
    \Lambda_H =
    d \cdot \parens*{ x_0 - (d-1) \frac{\chi_{x_0}(H)}{\chi'_{x_0}(H)} }
    \label{eq:lowerbound_dsq}
\end{equation}
where $\chi_\bullet(H)$ denotes the characteristic polynomial of the
adjacency matrix $A_H$ of $H$ evaluated at $x_0$, which is the
largest root to the equation
\begin{equation}
    (d-1) \cdot \chi_x(H) \cdot \chi''_x(H) = (d-2) \cdot \chi'_x(H)^2
    \label{eq:x0_def}
\end{equation}
cf.~\cite[Theorem~2.1]{CohenCMP25}. The authors also conjecture that
the lower bound from \eqref{eq:lowerbound_dsq} is saturated
by random graphs
drawn from the matching model except with vanishing probability. Namely,
\begin{equation}
    \Pp_{\rv M^N}\sqbr*{ \lambda\parens*{ \Grand\parens*{\rv M^N} \dsq
    H } \ge \Lambda_H + \ee } = o_N(1)
    \label{eq:dsq_conj}
\end{equation}
where $\Grand(\rv M^N)$ denotes a $d$-regular graph drawn from the
matching model, cf.~\cite[Conjecture~2.6]{CohenCMP25}. That is to
say, writing $\Arand(\rv M^N)$ for the adjacency matrix of
$\Grand(\rv M^N)$, we write
\[ \Arand(\rv M^N) = \sum_{i= 1}^d M_i^N\]
where $\rv M^N$ is a tuple of $d$ iid uniformly random matching
matrices. In this section, we prove the conjecture stated in
\eqref{eq:dsq_conj}. Namely, we will show
\begin{theorem}\label{thm:dsq_conj}
    For every $k$-regular vertex-transitive graph $H = ([d], E)$
    (where $d \ge 3, k \ge 1)$ and
    even $N \in \NN$, we have
    \[ \Pp_{\rv M^N}\sqbr*{ \lambda\parens*{ \Grand\parens*{\rv M^N}
    \dsq H } \ge \Lambda_H + \ee} = o_N(1) \]
    where $\rv M^N = (M_1^N, \ldots, M_d^N)$ denotes a $d$-tuple of
    iid uniformly random matching matrices and $\Lambda_H$ is as
    defined in \eqref{eq:lowerbound_dsq}.
\end{theorem}
To establish \cref{thm:dsq_conj}, we will first show that that the
adjacency matrix of the
derandomized square $G \dsq H$ can be thought as a matrix polynomial
in permutation (or matching) matrices --- this will follow from a
close inspection of the analysis of the derandomized square
construction occurring in \cite{RozenmanV05}. We will carry out this
task in \cref{sub:derandomized_square_as_a_matrix_polynomial}.

This already implies
that the spectral expansion of the derandomized square of a random
matrix from the matching model with $H$ can be analyzed using norm
considerations in the free-product $\Gam_d = \ZZ_2^{*d}$. It turns out
the analysis of \cite{CohenCMP25} is already sufficient to bound the
norm of the limiting object in $\Cred(\Gam_d)$. We will recall their
result and present the proof of \cref{thm:dsq_conj} in
\cref{sub:spectral_bounds_for_the_derandomized_square_of_the_infinite_tree}.

As a consequence of the norm bound, using
\cref{thm:constructive_convg} we also have an efficient algorithm for
constructing a family of $d$-regular near-Ramanujan graphs for which
we have $\lambda(G \dsq H) \le \Lambda_H + \ee$ for graphs $H$
satisfying the consequence of \cref{thm:dsq_conj}.
\begin{corollary}\label{cor:dsq_constr}
    For any given $N, d, k \in \NN$ ($d \ge 3, k \ge 1$) and $\ee >
    0$, there is a $d$-regular graph $G_{N'} = (V_{N'},
    E_{N'})$ on $N' \in [N, (1+ o_N(1))N]$ vertices such that
    \begin{enumerate}[(i)]
        \item $\lambda(G_{N'}) \le 2\sqrt{d-1} + \ee$, i.e., $G_{N'}$
            is almost Ramanujan.
        \item $\lambda(G_{N'} \dsq H) \le \Lambda_H + \ee$
            simultaneously for every $k$-regular vertex transitive
            graph $H = ([d], E')$.
    \end{enumerate}
    Moreover, the graph $G_{N'}$ can be constructed in time $N^{O_{d,
    k, \ee}(1)}$.
\end{corollary}
We note that \cref{cor:dsq_constr} is closely related to
\cite[Theorem 6.2]{CohenCMP25}, which shows that for every $d, N \in
\NN$, in the setting of \cref{cor:dsq_constr} there are
$d$-regular $H$-local graphs $G_{Nd}$ on $Nd$
vertices such that $\lambda(G_{Nd}) \le \Lambda_H$. The class of
$H$-local graphs (\cite[p.~35]{CohenCMP25}) contains the derandomized
square $G \dsq H$ for some base graph $G$, but they need not always
be a derandomized square. Our result trades the exact bound
$\Lambda_H$ for an approximation $\Lambda_H + \ee$, but ensures that
such graphs can be constructed as the derandomized square using the
same base graph $G_{N'}$ simultaneously for all vertex-transitive
graphs $H$. In particular, for these graphs the derandomized square always
performs like the average case (\cref{thm:dsq_conj}) as opposed to
the worst case bound \eqref{eq:dsq_expansion}.

\subsection{Derandomized Square as a Matrix
Polynomial}\label{sub:derandomized_square_as_a_matrix_polynomial} 
In \cite{RozenmanV05},
\eqref{eq:dsq_expansion} is proven by decomposing the random walk
matrix $M_S = A_S/(dk)$ of $S = G \dsq H$ by writing
\begin{equation}
    M_S = P \widetilde A \widetilde B \widetilde A
    L\label{eq:dsq_decomposition}
\end{equation}
where
\begin{enumerate}[(i)]
    \item $\widetilde A \in \RR^{Nd \times Nd}$ is defined by
        $\widetilde A_{(u, a); (u', a')} = \One\set*{a = a' \land
        \Out_G(u, a) = u'}$ for all $u, u' \in [n], a, a' \in [k]$,
    \item $\widetilde B \in \RR^{Nd \times Nd}$ satisfies $\widetilde
        B = I \tensor M_H$ where $M_H$ denotes the random
        walk matrix
        of the auxiliary graph $H$,
    \item $P \in \RR^{N \times Nd}$ and $L \in \RR^{Nd \times N}$ are
        projection and lift operators respectively, defined by
        \[ (P z)_u = \sum_{a = 1}^d z_{u, a}\qquad\textrm{and}\qquad
        L f = f \tensor \frac{\One_{[d]}}{d}  \]
\end{enumerate}
cf.~\cite[Proof of Theorem~4.4]{RozenmanV05}. We note that writing $A_G$
as in \eqref{eq:a_decomp}, it is easy to verify that
\[ \widetilde A = \sum_{j = 1}^d \lambda(\pi_j) \tensor
\One_j \One_j^\top\]
where $\One_j \One_j^\top$ is the matrix which equals 1 on the
$(j,j)$-th element and 0 elsewhere. Consequently,
\begin{equation}
    \widetilde A \widetilde B \widetilde A = \sum_{i, j=  1}^d
    \lambda(\pi_i) \lambda(\pi_j) \tensor \frac{\One\set*{\set*{i,j} \in
    E(H)}}{k} \cdot \One_i \One_j^\top.
    \label{eq:aba_decomp}
\end{equation}
Observing that $P$ and $L$ only act on the second tensor component in
\eqref{eq:aba_decomp}, we can compute
\begin{lemma}\label{lem:compressed_form}
    Let $S = G \dsq H$, and $\pi_1, \ldots, \pi_d$ be the matrices
    occurring in the decomposition \eqref{eq:a_decomp}. Then,
    \begin{equation}
        M_S = \frac{1}{dk} \sum_{\set*{i,j} \in E(H)} \parens*{
        \lambda(\pi_i) \lambda(\pi_j) + \lambda(\pi_j) \lambda(\pi_i)}.
        \label{eq:compressed_dsq}
    \end{equation}
    In particular, writing $\rv \Pi^N = (\lambda(\pi_1), \ldots,
    \lambda(\pi_n))$ for the matrix representation of the
    permutations occuring in \eqref{eq:a_decomp}, we have that
    \begin{equation}\label{eq:dsq_polynomial} A_S = Q_H(\rv \Pi^N) =
        \sum_{\set*{i,j} \in E(H)}
        \parens*{ \Pi_i^N \Pi_j^N + \Pi_j^N \Pi_i^N}.
    \end{equation}
    In particular, the random walk matrix $M_S$ can be expressed as a
    $*$-polynomial of degree 2 in $\rv \Pi^N$.
\end{lemma}
\begin{remark}\label{rem:switcheroo}
    When $G := \Grand(\rv M^N)$ is drawn from matching model, each
    $\pi_j$ occurring  in
    \cref{lem:compressed_form} is an involution, thus $\pi_j^{-1} =
    \pi_j$ for all $j \in [d]$. In particular, this means that
    $\Pi_j^N = \Pi_j^{N*}$, and thus the polynomial $Q_H(\rv \Pi^N)$
    from \eqref{eq:dsq_polynomial} can be thought as a self-adjoint
    $*$-polynomial by choosing the representation
    \[ Q_H(\rv \Pi^N) = \sum_{\set*{i,j} \in E(H)}
    \parens*{\Pi_i^N \Pi_j^N + \Pi_j^{N*} \Pi_i^{N^*}}.\]
\end{remark}
\begin{remark}
    In our simplified formula we have assumed that the graph has
    no-loops, in the presence of loops we can work with $Q_H(\rv
    \Pi^N) = \sum_{i,j = 1}^d (A_H)_{i,j} \cdot \Pi^N_i \Pi^{N*}_j$.
    Allowing $H$ to have loops, would force us to assume also that
    every vertex in $H$ has a non-loop edge to ensure that the
    results advertised in the beginning of the present section
    continue to hold.
\end{remark}

\subsection{Spectral Bounds for the Derandomized Square of the
    Infinite
Tree}\label{sub:spectral_bounds_for_the_derandomized_square_of_the_infinite_tree}
Let now $m_1,\ldots, m_d$ be the canonical involutive generators of
$\Gam_d := \ZZ_2^{*d}$
(the free-product of $\ZZ_2$ with itself, $d$-times). Throughout this
section we will assume $d \ge 3$ and the degree $k$ of $H$ is $k\ge
1$. We define the
polynomial the non-commutative $*$-polynomial $Q_H \in \CC^*\langle
u_1,\ldots, u_d\rangle$ of degree 2 by
\begin{equation}
    Q_H(u_1,\ldots, u_d) =  \sum_{\set*{i,j} \in E(H)}
    \parens*{ u_i u_j + u_j^* u_i^*}\qquad\textrm{where}\qquad u_i =
    \lambda(m_i)\quad\textrm{for all}\quad i =1,\ldots, d.
    \label{eq:quu_def}
\end{equation}
As before we will write $\rv u = (u_1,\ldots, u_d)$. Notice that
since $m_i^{-1} = m_i$, we have $u_i^* = u_i$, and as
such the polynomial $Q_H$ defined in \eqref{eq:quu_def} is
self-adjoint. Since the infinite tree $t$ can be equipped with the
adjacency operator $t = (u_1 + \cdots + u_d)$ in light of
\cref{lem:compressed_form,rem:switcheroo} $Q_H(\rv u)$ can be
identified with
the adjacency operator of the derandomized square $\Tre_d \dsq H$. To
this end, we will rely on \cref{lem:norm_def} to compute the norm
$\norm*{Q_H(\rv u)}_{\Cred(\Gam_d)}$.
\begin{proposition}\label{prop:analy_comb}
    Let $Q_H(\rv u)$ be as defined in \eqref{eq:quu_def}. Then,
    writing $ Q_H(\rv u) \in \Cred(\Gam_d)$, we have
    \begin{align*} \norm*{s}_{\Cred(\Gam_d)} &~=~ \lim_{p \to \infty}
        \parens*{\tau\parens*{ s^{2p} }}^{\frac{1}{2p}},\\
        &~=~ \lim_{p \to
        \infty} \Big[ \textsc{Number of Closed Walks in $\Tre_d \dsq H$
        of Length $2p$ Originating at $e$} \Big]^{\frac{1}{2p}}
    \end{align*}
    where $e \in \Gam_d$ is the identity word.
\end{proposition}
\begin{proof}
    By \cref{lem:norm_def}, we know
    \[ \norm*{s}_{\Cred(\Gam_d)} = \lim_{p \to \infty}
    \parens*{\tau\parens*{ (s^* s)^p }}^{\frac{1}{2p}}.\]
    The first equality thus follows by the self-adjointness of $s$,
    i.e.~$s^* = s$, which is evident from \eqref{eq:quu_def}. The
    second equality follows from the definition of the canonical
    trace (q.v.~\cref{xmpl:group_algebras}). We have, by identifying
    the root of $\Tre_d$ with the identity word $e \in \Gam_d$,
    \begin{align*} \tau\parens*{s^{2p}} &~=~ \inpr*{s^{2p} \delta_e,
        \delta_e},\\
        &~=~\Big[\textsc{Number of Closed Walks in $\Tre_d \dsq H$ of
                Length $2p$ Originating at
        $e$}\Big]
    \end{align*}
    which follows since $s$ is the adjacency operator of $\Tre_d \dsq
    H$, per our previous discussion.
\end{proof}
\cref{prop:analy_comb} reduces the problem of bounding
$\norm*{s}_{\Cred(\Gam_d)}$ to the problem of counting the number of
closed walks in $\Tre_d \dsq H$. Now, we describe the analysis of
\cite{CohenCMP25} that led to their lower bound
\eqref{eq:lowerbound_dsq}.

First, we give a description of the
infinite graph $\Tre_d \dsq H$ whose adjacency operator is given by
$s = Q_H(\rv u)$. We root the infinite tree $\Tre_d$ at $e
\in \Gam_d$ (identity word), and for every vertex $v \in \Gam_d$, we
label its $d$ children $\set*{v_1,\ldots, v_d}$ consistently with
labels in $[d]$,
i.e.~if $w = v_i$ then $w_i = v$ for each $v, w \in \Gam_d$. Among
the neighbors of each $v \in \Gam_d$ in $\Tre_d$, we place a graph
$H$ which retains the adjacency structure of $H$, i.e.~$\set*{v_i,
v_j}$ is an edge we add if $\set*{i, j} \in H$.  Notice that the
edges $\set*{v_i, v_j}$ are precisely the set of edges of $\Tre_d
\dsq H$. Notice that edges in $\Tre_d \dsq H$ can be identified with
a (subset) of the walks in
in $\Tre_d$ of length 2.

Equivalently, we can build an infinite graph, by starting at the root and
attaching it $d$-copies of $H$, at each newly attached vertex we can add
$(d-1)$-copies of $H$ and recurse by attaching $(d-1)$-copies to all
attached vertices also. In particular, we can say that the root has
$d$ forward copies of $H$ attached to it, and the remaining vertices
has one parent copy attached to them (via their parent in $\Tre_d$)
and $(d-1)$-forward copies (via their children in $\Tre_d$).

We will write $\mathbb X$ for the connected component of $\Tre_d
\dsq H$, after
removing one of the copies adjacent to the root, i.e., every vertex
(including the root) in $v \in V(\mathbb X)$ are contained in
$(d-1)$-forward copies of $H$, and non-root vertices are contained in
a parent copy of $H$. This is the pruned graph that
\citeauthor{CohenCMP25} \cite{CohenCMP25} consider. As this graph
contains arbitrarily large neighborhoods agreeing with $\Tre_d \dsq
H$, we can the asymptotic estimates for cycle
counts in $\mathbb X$ to bound the corresponding cycle counts in
$\Tre_d \dsq H$.

The idea of \cite{CohenCMP25} is bounding coefficients of the
generating function
\begin{equation}\label{eq:c_def} C_H(z) = 1 + \sum_{p \ge 1} c_p z^p
\end{equation}
where $c_p$ denotes the number of length $p$ closed walks originating
at $e$ in $\mathbb X$.

Concretely, using methods of analytic combinatorics (in the spirit of
\cite{FlajoletS09}) they prove
\begin{theorem}[\cite{CohenCMP25}, Corollary 5.7]\footnote{Using the
        notation of
        \cite{CohenCMP25}, we have $\zeta_H(z)=zC_H(z)$, and hence
        $c_{2p}$ is the coefficient of $z^{2p+1}$ in $\zeta_H(z)$.
        When $H$ is bipartite, even powers in $\zeta_H$ vanish,
        so $c_{2p}$ are precisely the non-zero coefficients in
        $\zeta_H$; in the
        non-bipartite case, the same estimate applies to even powers as
    well.}\label{thm:cycle_pruned}
    Let $C(z)$ be as above. Then,
    \[ \lim_{p \to \infty} \parens*{ c_{2p}}^{\frac{1}{2p}}
    = \Lambda_H\]
    where $\Lambda_H$ is as defined in \eqref{eq:lowerbound_dsq}.
\end{theorem}
\cref{thm:cycle_pruned} bounds the growth rate of $c_{2p}$, the
number of closed walks of length $2p$ in the pruned infinite graph
$\mathbb{X}$ originating at $e$. It turns out, we can relate this
number to the number of closed walks in $\Tre_d \dsq H$ which
originate at $e$ and have length $2p$, which allows us to bound
$\norm*{s}_{\Cred(\Gam_d)}$ by appealing to \cref{prop:analy_comb}.

\begin{corollary}\label{cor:s_norm}
    Let $s = Q_H(\rv u)$ where $Q_H(\bullet)$ and $\rv u$ are as
    defined in \eqref{eq:quu_def}. Then,
    \[ \norm*{s}_{\Cred(\Gam_d)} = \Lambda_H.\]
\end{corollary}
\begin{proof}
    Notice that every path emerging from
    $v$ in $\mathbb{X}$ is also a path emerging from $v$ in the graph
    $\Tre_d \dsq H$, since $\mathbb{X}$ is a subgraph of $\Tre_d \dsq
    H$. Consequently, we have by \cref{thm:cycle_pruned}
    \begin{equation}\label{eq:simple_ineq} c_{2p} \le
        \tau\parens*{s^{2p}} \Longrightarrow \lim_{p \to
        \infty} c_{2p}^{\frac{1}{2p}} = \lambda_H \le
        \norm*{s}_{\Cred(\Gam_d)}
    \end{equation}
    as by \cref{prop:analy_comb} we know that $\tau(s^{2p})$ counts
    closed walks in $\Tre_d \dsq H$ of length $2p$ of length $2p$ and
    $c_{2p}$ counts such walks in a subgraph.

    Now, we fix a vertex $v \in \Gam_d$ such that every path emerging
    from $v$ in the graph $\Tre_d \dsq H$ of length $2p$ is also
    present in the
    graph $\mathbb{X}$.\footnote{Each edge of $\mathbb X$ corresponds
        to a walk of
        length two in $\Tre_d$. Starting at the root, we can repeatedly
        traverse a non-loop edge in a copy of $H$ different from the
        copy we came from. To every vertex $d-1$ forward copies
        are attached, so we can reach vertices arbitrarily far
        removed from the root.  Vertex-transitivity ensures that
        there is always a forwards edge
        we can traverse, so long as there is at least one non-loop
    attached to every vertex.} Since both $H$ and $\Tre_d$ are vertex
    transitive, the number of closed walks of length $2p$ from $v$ in
    $\Tre_d \dsq H$ is the same as the number of closed walks of
    length $2p$ emerging $e$. Moreover, by our choice of $v$, this
    equals the number of closed walks in $\mathbb X$ emerging from $v$.

    This enables us thus to upper bound $\tau(s^{2p})$ by closed
    walks in $\mathbb X$: Writing $\Delta$ for the distance from $e$
    to $v$ in $\mathbb X$, we can always move from $e$ to $v$
    following along this path, follow any $q$ (potentially distinct)
    closed walks emerging from $v$ of length $2p$, and follow
    the path $v$ to $e$ in reverse direction. The number of such
    paths is $c_{2pq + 2\Delta}$ and by the preceding discussion
    $c_{2pq + 2\Delta}$ can be used to upper bound
    $\tau(s^{2p})^q$. Formally,
    \[ \tau\parens*{s^{2p}}^q \le c_{2pq + 2\Delta} \Longrightarrow
        \tau\parens*{s^{2p}}^{\frac{1}{2p}} \le \parens*{c_{2pq +
    2\Delta}}^{\frac{1}{2pq}}.\]
    For fixed $p$ and $\Delta$, letting $q \to \infty$, we get by
    \cref{thm:cycle_pruned}
    \[ \lim_{q \to \infty} \parens*{c_{2pq +
    2\Delta}}^{\frac{1}{2pq}} = \Lambda_H.  \]
    Letting $p \to \infty$ (and therefore $\Delta \to \infty$ also)
    also, we obtain
    \[ \norm*{s}_{\Cred(\Gam_d)} = \lim_{p \to \infty}
    \tau\parens*{s^{2p}}^{\frac{1}{2p}} \le \Lambda_H\]
    which implies the statement together with \eqref{eq:simple_ineq}.
\end{proof}
Proof of \cref{thm:dsq_conj} is now immediate from
\cref{thm:asymptotic_bc} by
the norm bound on
$s$ given in \cref{cor:s_norm} and
from the fact that $s := Q_H(\rv u)$ (\cref{lem:compressed_form})
and $Q_H$ is a quadratic $*$-polynomial.

Proof \cref{cor:dsq_constr} will require the additional observation
that the adjacency of the operator of $\Tr_d$ can be identified with $t
= u_1 + \cdots + u_d$ where $u_i$ are as in \eqref{eq:quu_def} and
$\norm*{t}_{\Cred(\Gam_d)} = 2\sqrt{d-1}$. In particular, we can
apply \cref{thm:constructive_convg} with the choice of polynomials
\[ \mathcal P = \set*{ Q_H\parens*{\rv u} \in \CC^*\langle \rv u \rangle :
    H~\textrm{is a $k$-regular vertex transitive graph on $[d]$} } \sqcup
\set*{T(\rv u) = \sum_{i = 1}^d  u_i} \]
where $Q_H$ and $\rv u$ as in \eqref{eq:quu_def}.
This is a collection of polynomials of maximum degree at most 2,
whose largest coefficient is bounded by 1 in absolute value.

\section{Permutation Conjugated Adjacency
Matrices}\label{sec:permutation_conjugated_adjacency_matrices} 
Let $\rv U^N = (U^N_1,\ldots, U_r^N)$ be a sequence iid uniform
permutation matrices and let $G_{N; 1},\ldots, G_{N; r}$ be a
sequence $d_i$-regular graphs on the same vertex set $V := [N]$.
Throughout this section, we will be concerned with joint moments of
permutation conjugated adjacency matrices,
\begin{equation}\label{eq:bee_def} B_i^N = U^N_i A_i^N
    U^{N*}_i\qquad\forall~i = 1,\ldots,r
\end{equation}
where $A_i^N$ denotes the adjacency matrix of the graph $G_{N; i}$
for each $i \in [r]$.

We will show a weak convergence theorem for the sequence $\rv B^N$: Let
\[ \Gam := \Gam_{\vec d, r} = *_{i = 1}^r *_{j = 1}^{d_i}\ZZ_2\qquad
    \textrm{and for all $i \in [r]$}\qquad t_i := \sum_{j =
1}^{d_i} \lambda(m_{i,j}) \]
where $m_{i,j}$ are the canonical free involutive elements which
generate $\Gam$. We will also write $\rv t=  (t_1,\ldots, t_r)$.

We note that this model of conjugated graphs we consider is general
enough to capture the rotating walks from
\cref{sec:rotating_infinite_trees_and_expanders}. Namely, in
\cref{sub:rotating_walks_reduction} we will show
\begin{proposition}\label{prop:rotating_redux}
    Let $G_N$ be any given $d$-regular graph on $N$ vertices with
    adjacency matrix $A^N$. Let $\rv U^N =
    (U^N_1,\ldots, U_{k}^N)$ be a sequence of iid uniform permutation
    matrices. Defining, $\rv B^N = (U_1^N A^N U_1^{N*}, \ldots, U_{k-1}^N
    A^N U_{k-1}^{N^*})$ we have that the spectral laws of the
    random matrices
    $\Arot\parens*{\rv U_{[2, k-1]}^N}$
    and $(B_1^N
    \cdots B_{k-1}^N) \cdot (B_1^N
    \cdots B_{k-1}^N)^*$ distribute the same.
    here $\Arot\parens*{\rv U_{[2, k-1]}^N}$ is the rotating walk
    as defined by
    \eqref{eq:rot_walk_def} and $\rv U^N_{[2, k-1]} = (U_2^N, \ldots,
    U_{k-1}^N)$.
\end{proposition}
%

For words $w = (i_1,\ldots, i_q) \in [r]^q$ we will use the notation
\[ w(\rv B^N) = B_{i_1}^N \cdots B_{i_q}^N\quad\tand\quad w(\rv t) =
t_{i_1}\cdots t_{i_q}.  \]
We will show that $\rv B^N$ weakly converges to $\rv t$
when the $d_i$-s are fixed and
(q.v.~\cref{sub:weak_and_strong_convergence}) and the graph $G_{N;
i}$ have girth $\omega(1)$ for all $i = 1,\ldots, r$. Concretely,
we will prove
\begin{theorem}\label{thm:weak_convergence}
    Suppose the graphs $G_{N; i}$ come from a family of $d_i$-regular
    graphs where
    $\girth(G_{N; i}) = \omega(1)$ and $d_i = O(1)$ for each $i =
    1,\ldots, r$. Then,
    $\rv B^N$ weakly converges to $\rv t$, where both $\rv B^N$ and
    $\rv t$ are defined as in \eqref{eq:bee_def}.

    In particular, for any word $w = (i_1,\ldots, i_q) \in [r]^q$,
    writing $q_j$ the number of times $j \in [r]$ occurs in $w$, we have
    \begin{equation}\label{eq:trace_conc}\Pp\sqbr*{ \Abs*{ \tr w(\rv
            B^N) - \tau(w(\rv t)) } \ge
            \prod_{i=1}^{r} d_i^{q_i} \cdot \parens*{ \frac{3q^2}{N} +
        \eta}} \le 2 \exp\parens*{-
        \frac{\eta^2 \cdot N}{32 q^2}}\qquad\forall~\eta > 0
    \end{equation}
    so long as $q_i < \girth(G_{N; i})$ for all $i = 1,\ldots, r$.
    Concretely, the choice of $\eta = 8q \sqrt{\frac{\log N}N}$
    ensures that the RHS is bounded by $2/N^2$ while ensuring that
    the deviation measured by \eqref{eq:trace_conc} tends to 0.
\end{theorem}
The proof of \cref{thm:weak_convergence} will rely on a
pointwise convergence argument for the expected value of $\Exp w(\rv
B^N)$ coupled with a standard concentration
argument for functions of independent permutations which can be
obtained via the modified log-Sobolev inequality for the random
transposition chain on permutations, \cite{Goel04}. We will present
this proof in \cref{sub:proof_of_weak_convergence}.

Concretely, we will make use of the following convergence result for
the expectation:
\begin{theorem}\label{thm:pointwise_convg}
    Let $w = (i_1,\ldots, i_q) \in [r]^q$ be given for $q \ge 1$.
    Then, with $\rv B^N = (B^N_1,\ldots, B_r^N)$ as defined in
    \eqref{eq:bee_def} we have
    \[ \Abs*{ \Exp \tr w(\rv B^N) - \tau(w(\rv t)) } \le \prod_{i =
    1}^r d_i^{q_i} \cdot \min\set*{1, \frac{3q^2}N }\]
    where for all $j \in [r]$, we have $q_j = \Abs*{\set*{a \in [q] :
    i_a = j}}$ and $q_j < \girth(G_{N; i})$.
\end{theorem}

We will prove \cref{thm:pointwise_convg} in
\cref{sub:pointwise_convergence} by using a probabilistic
argument. First, we observe that we can think of $\tr w(\rv B^N)$ and
$\tau(w(\rv t))$ as the probabilities that a random walk reaches its
initial vertex after $q$ steps.

Consider an edge colored graph $\hat G$ whose set of vertices
coincides with $V$ and the set of edges correspond to the edges in
$E_i$ (of $G_{N; i})$ relabeled according to the random permutations.
Each edge coming from a relabeling of $E_i$ is identified with the
color $i$. Then, we observe
\[ \frac{\Exp\tr w(\rv B^N)}{\prod_{i = 1}^r d^{q_i}}  =
    \Pp\sqbr*{\textrm{a colored walk according to $w$ in $\hat G$ returns
to it's initial vertex}}\]
where the probability is taken both over the choice of random
relabelings and the walk. Here, we mean that a colored walk according
to $w = (i_1,\ldots, i_q)$ takes an edge with the color $i_j$ at the
$j$-th step for all $j \in [q]$.

Similarly, we can think of $\tau(w(\rv t))$ as being analogously
defined, but in the infinite tree $\hat{\Tre} = \Tre_{D}$ where $D =
\sum_{i=1}^r d_i$ and every vertex has $d_i$ neighbors with the
color $i$. Thus,
\[\frac{\tau(\rv t)}{\prod_{i = 1}^r d^{q_i}} =
    \Pp\sqbr*{\textrm{a colored walk starting from the root according
            to $w$ in $\hat \Tre$ returns
to it's initial vertex}}.\]
It turns out, the difference between these probabilities can
conveniently be bounded by analyzing the \enquote{shapes} which
underlie these walks.

A known consequence of weak convergence is an asymptotic lower
bound on $\norm*{
Q(\rv B^N) }$ for non-commutative $*$-polynomials $Q \in \CC^*\langle
x_1,\ldots, x_r\rangle$ by $\norm{Q(\rv t)}_{\Cred(\Gam)}$ which
holds with probability $1- o(1)$, q.v.~\cref{lem:weak_bd}.
\begin{corollary}\label{cor:weak_bd}
    Suppose the graphs $G_{N; i}$ come from a family of $d_i$-regular
    graphs where
    $\girth(G_{N; i}) = \omega(1)$ and $d_i = O(1)$ for each $i =
    1,\ldots, r$. Then, for any non-commutative $*$-polynomial $Q \in
    \CC^*\langle x_1,\ldots, x_r\rangle$ we have
    \[ \norm*{Q(\rv B^N)} \ge \norm*{Q(\rv t)}_{\Cred(\Gamma)} -
    o(1)\qquad\textrm{with probability at least}\qquad 1 - o(1)\]
    where $\rv B^N$ and $\rv t$ are as defined in \eqref{eq:quu_def}.

    The same consequence holds with $\rv B^N$ replaced by $\rv C^N =
    (\left. B_1^N\right|_{\perp 1}, \ldots, \left. B_r^N
    \right|_{\perp 1})$.
\end{corollary}
\begin{proof}
    The statement for $Q(\rv B^N)$ is an immediate consequence of
    \cref{lem:weak_bd}. The statement for $\rv C^N$ follows from the
    observation, that for any word $w = (i_1,\ldots, i_q)$ we have
    \begin{equation}\label{eq:b-c-close} \tr w(\rv C^N) = \frac{\Tr
            w(\rv B^N) - \prod_{i=1}^r
        d_i^{q_i}}{N} = \tr w(\rv B^N) - \frac{\prod_{i=1}^r d_i^{q_i}}{N}
    \end{equation}
    since $1$ is an eigenvector of
    \[ w(\rv B^N) = B_{i_1}^N \cdots B_{i_q}^N\]
    with eigenvalue $\prod_{i=1}^r d_i^{q_i}$. Thus, when $q$,
    $d_i$-s and $c_i$-s are constant \eqref{eq:b-c-close} verifies
    that $\tr w(\rv B^N)$ and $\tr w(\rv C^N)$ are at most $o(1)$
    apart. Consequently, the weak convergence result for $\rv B^N$
    extends to $\rv C^N$ also.
\end{proof}

It turns out when the polynomials $Q$ we consider have real non-negative
coefficients, we can also prove a lower bound for $\norm*{\left. Q(\rv
B^N)\right|_{\perp 1}}$ or which holds with probability 1.

\begin{theorem}\label{thm:general_lowerbound}
    Let $Q \in \RRp^*\langle x_1,\ldots, x_r\rangle$ be a
    non-commutative $*$-polynomial with non-negative real
    coefficients. Then, for every $N \ge 2$, we have
    \[ \norm*{Q(\left.\rv B^N)\right|_{\perp 1}} \ge \norm*{Q(\rv
        t)}_{\Cred(\Gamma)} -
    o_N(1) \]
    with probability 1, where $\rv B^N$ and $\rv t$ are as defined in
    \eqref{eq:quu_def}
    and we assume $d_i=  O(1)$ for all $i = 1,\ldots, r$.
\end{theorem}
We will prove this result in
\cref{sub:spectral_lower_bounds_for_permutation_conjugated_matrices}.

In particular, in conjunction with
\cref{prop:rotating_redux,thm:inf_rot_tree} the
above result allows us to recover the lower bound \eqref{eq:rot_lbd}
of \cite{CohenCMP25}.
\begin{corollary}
    Let $\rv \Pi^N = (\Pi_1^N, \ldots, \Pi_{k-1}^N)$ be any sequence of
    permutation matrices and $G_N$ any $d = O(1)$-regular graph on $N$
    vertices with adjacency matrix $N$, then
    \[ \lambda\parens*{\Arot(\rv \Pi^N)} \ge     \parens*{\frac{k+1}{k}}^k
    \cdot (k+1) \cdot (d-1) d^{k-1} - o_N(1).\]
\end{corollary}
\subsection{Convergence in Expectation for
Monomials}\label{sub:pointwise_convergence} 
Recall that we have defined the edge colored (random) graph $\hat G$
on the vertex set $V$ with colored edges $E_i$ (with the color $i$)
each corresponding to a (random) relabeling of the original set of
edges in $G_{N; i}$.
We have observed that,
\begin{equation}\label{eq:hatg_formula} \frac{\Exp\tr w(\rv
    B^N)}{\prod_{i = 1}^r d^{q_i}}  =
    \Pp\sqbr*{\textrm{a colored walk according to $w$ in $\hat G$ returns
    to it's initial vertex}}
\end{equation}
where the probability is taken both over the choice of random
relabelings and the walk. We recall that a colored walk according
to $w = (i_1,\ldots, i_q)$ takes an edge with the color $i_j$ at the
$j$-th step for all $j \in [q]$.

Similarly, we argued that $\tau(w(\rv t))$ measuring a similar probability
defined, but in the infinite tree $\hat{\Tre} = \Tre_{D}$ where $D =
\sum_{i=1}^r d_i$ and every vertex has $d_i$ neighbors with the
color $i$. Thus,
\begin{equation}\label{eq:hattree_formula}\frac{\tau(\rv t)}{\prod_{i
    = 1}^r d^{q_i}} =
    \Pp\sqbr*{\textrm{a colored walk starting from the root according
            to $w$ in $\hat \Tre$ returns
    to it's initial vertex}}.
\end{equation}
Instead of counting such walks directly, we will take a probabilistic
coarsening. For each $i = 1,\ldots, r$ we will call the shape of
a (not-necessarily closed) walk
$(v_0,\ldots, v_q)$ the
sequence obtained by relabeling
by their order of appearance.
So, $(5, 5, 3, 1, 7, 5)$ has the shape $(0, 0, 1, 2, 3, 0)$. Then,
notice that the shape of any walk in the tree $\hat \Tre$ will
correspond to a tree after taking the union of all the visited edges
(while keeping their color) and ignoring edges traversed multiple
times. It turns out, this insight is enough to count the number of
(colored, but not closed) walks corresponding to a shape. Since the
shape is enough to identify whether a walk is closed or not, we can
reduce the problem of bounding the difference between the
probabilities of walks returning to origin vertex simply by bounding
the total variation distance between the distributions of shapes
which occur in walks of $\hat G$ and $\hat \Tre$ which are weighted
by the number of walks which gave rise to them.

We recall that when $\mu$ and $\nu$ are supported on some joint set
$\mathcal X$, i.e., $\supp(\mu), \supp(\nu) \subset \mathcal X$ we have
\begin{equation}\label{eq:dtv_def}
    \norm*{\mu - \nu}_{\mathrm{tv}} = \frac{1}2 \sum_{x \in \mathcal
    X} |\mu(x) - \nu(x)| = \max_{S \subset \mathcal
    X}\Abs*{ \mu(S) - \nu(S) }
\end{equation}
We will need the following convenient formula,
\begin{fact}[\cite{LevinPW17}, Eq.~4.13]\label{fac:tvb}
    For any $\mu, \nu$ and finite collection $\mathcal X$, such that
    $\supp(\mu), \supp(\nu) \subset \mathcal X$, we have that
    \[\norm*{\mu - \nu}_{\mathrm{tv}} = 1 - \sum_{x \in \mathcal X}
    \min\set*{\mu(x), \nu(x)}.\]
\end{fact}
Equipped with \cref{fac:tvb}, we will show that the distribution of
shapes over the tree $\hat \Tre$ and the random graph $\hat G$ are
very close to each other statistically.
\begin{lemma}\label{eq:stat-dist}
    Writing $\tilde \beta(\hat \Tre, \bullet)$ for the $\tilde
    \beta(\hat G, \bullet)$ for the distributions of the shapes
    weighted by the number of walks (not-necessarily closed) which
    give rise to them, we have
    \[ \norm*{ \tilde \beta(\hat \Tre, \bullet) - \tilde\beta(\hat G,
    \bullet)}_{\mathrm{tv}} \le     \frac{3q^2}{N}\]
    where we recall that $\hat G$ itself is a random graph and
    $\beta(\hat G,
    \bullet)$ also accounts for this randomness.
\end{lemma}
We note that \cref{eq:stat-dist} immediately implies
\cref{thm:pointwise_convg} by observing that the normalizing factor
in \eqref{eq:hatg_formula} and \eqref{eq:hattree_formula} is
$\prod_{i=1}^r d_i^{q_i}$.
\begin{proof}[Proof of \cref{eq:stat-dist}]
    Let $T$ be an arbitrary tree. Before we start, we fix some notation. For
    each $i \in [r]$, we
    will write $F_i$ for the forests contained in $T$ which are
    colored with the color $i$. (We can assume w.l.o.g that no such
        $F_i$ is empty and we will omit the $F_i$ which correspond to
    isolated vertices.)

    We will write $\ell$ for the number of vertices contained in $T$,
    $f_i$ for the number of vertices in each $F_i$, and $\kappa_i$
    for the number of connected components for each $i \in [r]$.
    We can verify,
    \begin{equation} \sum_{i=1}^r (f_i - \kappa_i) = \ell - 1;\quad
        \ell \le q + 1;
        \quad \sum_{i=1}^r f_i \le 2q;\quad \kappa_i \le q.\label{eq:simple}
    \end{equation}

    Any connected component (isomorphic to) $C$ of $F_i$ starting
    from the root can
    be traversed by
    \[ \alpha_{i, C} = (d_i)_{\deg_C(e)} \cdot \prod_{v \ne e} (d_i -
    1)_{\deg_C(v)- 1}\]
    walks emanating from $e$ of color $i$ where $(n)_k = n (n-1)
    \cdots (n-k + 1)$ is the falling factorial. This can be thought
    as the number of way doing a BFS on $\Tre_{d_i}$ starting
    from the root while ensuring that we preserve the degree counts
    of $C$. Since we the shape structure tells us when we do
    backwards edges, this precisely counts the number of walks in
    tree $\hat T$ with a given tree shape $T$.

    In particular, we can compute
    \begin{equation}\label{eq:prophmct} \tilde\beta(\hat \Tre, T) =
        \frac{\prod_{i = 1}^r \prod_{C \le
        F_i} \alpha_{i, C}}{\prod_{i=1}^r d_i^{q_i}}
    \end{equation}
    since the denominator is simply the count of walks in the tree of
    length $q$ emanating from the root of $\hat \Tre$.

    Writing $\mathrm{inj}(F_i, G_{N; i})$ for the number of ways
    vertices of the
    forest $F_i$ can be mapped invectively to $G_{N; i}$ while
    ensuring that edges are mapped to edges, we have
    \begin{equation} \tilde \beta(\hat G, T) = \frac{(N)_\ell}{N
            \cdot \prod_{i=1}^{r}
        d_i^{q_i}}\prod_{i = 1}^r \frac{\mathrm{inj}(F_i,
        G_{N; i})}{(N)_{f_i}}\label{eq:ghomct}
    \end{equation}
    where $(N)_\ell$ is the number of ways to map $\ell$ vertices
    injectively to $[N]$, $N \prod_{i= 1}^r d_i^{q_i}$ counts the
    total number of walks, $\mathrm{inj}(F_i, G_{N; i})/(N)_{f_i}$ is
    the probability that all required edges exist to have traversed
    $F_i$. Since different colors use independent permutations the
    formula follows.

    Next, we note that we can compare $\tilde \beta(\hat G, T)$ to
    $\tilde \beta(\hat \Tre, T)$ by means of the following
    experiment: We can try to embed a component $C$ of
    $F_i$ in $G_{N; i}$ by picking the location of any one designated root
    vertex $e\in C$ at random to some uniform $v$ and mapping it
    injectively to a tree by mapping the neighbors of $e_i$
    injectively to the neighbors of $e_i$ and continue this process
    recursively until we have exhausted every edge in $C$. We observe
    that this tree growing procedure ensures injectivity in each
    component since every component has at most $q_i < \girth(C_i)$
    edges, i.e., non-injectivity implies the existence of a
    non-backtracking cycle of length at most $q_i$. The root
    of each $e$ is placed randomly and independently and components
    (subtrees) are also grown randomly, thus the probability that any two
    arbitrary vertices in different components are mapped to the same
    vertex is exactly $1/N$. By a simple union bound argument, we
    therefore see that the probability that this produces a
    non-injective mapping preserving edges is at most
    $\binom{f_j}2/N.$ This implies,
    \[ \mathrm{inj}(F_i, G_{N; i}) \ge N^{\kappa_i}\prod_{C \le F_i}
    \alpha_i(C) \cdot \parens*{1 - \frac{f_j^2}{2N}}\]
    since $N^{\kappa_i} \prod_{C \le F_i} \alpha_i(C)$ is the total
    labelings which preserve the edges in each $C$ and are injective
    when restricted to each $C$, and
    $1- f_j^2/2N$ is a lower bound on a random labeling such this
    being globally injective. Thus, we get
    \begin{equation} \prod_{i= 1}^r \frac{\mathrm{inj}(F_i; G_{N;
        i})}{N^{\kappa_i}
        \prod_{C \le F_i}\alpha_i(C)} \ge 1 - \sum_{i =
        1}^s \frac{f_i^2}{2N}.\label{eq:magic}
    \end{equation}
    In particular by \cref{eq:prophmct,eq:ghomct} and the first item
    in \eqref{eq:simple},
    \begin{align*}
        \frac{\tilde\beta(\hat G, T)}{\tilde\beta(\hat\Tre, T)}
        &~=~\frac{(N)_\ell}{N^\ell} \prod_{i=1}^r
        \frac{N^{f_i}}{(N)_{f_i}} \cdot \frac{\mathrm{inj}(F_i; G_{N;
        i})}{N^{\kappa_i} \prod_{C \le F_i} \alpha_i(C)}.
    \end{align*}
    By using \eqref{eq:magic}
    \begin{align*}
        \frac{\tilde\beta(\hat G, T)}{\tilde\beta(\hat\Tre, T)}
        &~\ge~\frac{(N)_{\ell}}{N^\ell} \cdot \parens*{1 -
        \sum_{i=1}^r \frac{f_i^2}{2N}},\\
        &~\ge~\parens*{1 - \frac{\binom{\ell}2}N} \cdot \parens*{1 -
        \sum_{i=1}^r\frac{f_i^2}{2N}}.
    \end{align*}
    By \eqref{eq:simple} and the Cauchy Schwarz inequality and
    $\sum_{i =1}^r f_i \le 2q$ we have
    \[ \frac{1}{2N} \sum_{i=1}^r f_i^2  \le \frac{1}{2N} \cdot
        \parens*{\sum_{i=1}^r f_i}^2 \le
        \frac{2q^2}{N}\qquad\textrm{and}\qquad
    \frac{\binom{\ell}2}{N} \le \frac{q^2}N.\]
    Thus,
    \[\frac{\tilde\beta(\hat G, T)}{\tilde\beta(\hat\Tre, T)}
    \ge 1 - \frac{3q^2}N.\]
    Finally, for any $T$ we obtain
    \[ \frac{\tilde\beta(\hat G, T)}{\tilde\beta(\hat\Tre, T)} \ge 1
        - \frac{3q^2}N \Longleftrightarrow \tilde\beta(\hat G, T)
        \ge\parens*{1 - \frac{3q^2}N} \cdot \tilde\beta(\hat \Tre, T).
    \]
    Thus, by \cref{fac:tvb}, we obtain
    \[ \norm*{\tilde\beta(\hat G, T) - \tilde\beta(\hat \Tre,
        T)}_{\mathrm{tv}} = 1 - \sum_{T} \min\set*{\tilde\beta(\hat G,
        T), \tilde\beta(\hat \Tre, T)} \le 1 - \sum_T \parens*{1 -
    \frac{3q^2}N} \beta(\hat\Tre, T) = \frac{3q^2}N\]
    where we note that it is sufficient to restrict the sum occurring
    in \cref{fac:tvb} to trees, since $\beta(\hat \Tre, S) =0$ for
    any non-tree shape.
\end{proof}
\subsection{Proof of Weak
Convergence}\label{sub:proof_of_weak_convergence} 
We will use the following concentration inequality for functions on
permutations, obtained from the modified log-Sobolev inequality for
the transposition chain \cite{Goel04} and the Herbst argument.  For a single
permutation it is essentially given by \cite[Corollary 3.1, Theorem
5.1]{Goel04}.
For multiple permutations, it is implied by a standard argument using
tensorization of entropy (\cite[Theorem 4.10]{BoucheronLM16}).
\begin{lemma}\label{lem:perm_conc}
    Let $\rv \pi = (\pi_1,\ldots, \pi_r)$ be a sequence of iid
    uniform permutations from $\Sym_N$, where $N \ge 2$, and let $F :
    \Sym_N^r \to \RR$ be a given function. Suppose changing any
    component $\pi_i$ of $\rv \pi$ with a copy of $\pi_i$ where two
    elements are transposed, i.e. $(a\ b) \pi_i$, changes $F$ by at
    most $L_i \ge 0$. Then,
    \[ \Pp_{\rv \pi}\sqbr*{ \Abs*{F - \mathbb E F} \ge \eta} \le 2
    \exp\parens*{- \frac{\eta^2}{2 (N-1)\sum_{i=1}^r L_i^2}}.\]
\end{lemma}

\begin{proof}[Proof of \cref{thm:weak_convergence}]
    Write $\hat B^N_i = B_i^N/d_i$ for all $i = 1,\ldots, r$ and define
    \begin{equation}\label{eq:effdeff}         F(\rv \pi) :=\tr\parens*{
        w(\rv{\hat B}^N)} = \frac{ \tr
        w(\rv B^N)}{\prod_{i=1}^r d_i^{q_i}}.
    \end{equation}
    where $\rv \pi = (\pi_1,\ldots, \pi_r)$ and $\pi_i$ is the
    permutation which is used to relabel each graph $G_i$. We will
    use the notation $\rv \pi^{i, a \lra b}$ for the copy of $\rv
    \pi$ where $\pi_i$ is replaced by $(a\ b) \pi_i$.
    Notice that switching from $\rv \pi$ to $\rv \pi^{i, a \lra b}$ only
    affects the  occurrences of the $i$-th matrix $B_i^N$ in
    \eqref{eq:effdeff}.
    Concretely, writing $T_{a \lra b} := T_{a \lra b}^N$ for the
    transposition
    matrix, it replaces $\hat B_i^N$ with $T_{a \lra b} \hat B_i^N
    T_{a \lra b}^*$ in \eqref{eq:effdeff}. Notice that,
    \[ \Delta_{i, a, b} := {  T_{a \lra b} \hat B_i^N
        T_{a \lra b} - \hat B_i^N } = (T_{a \lra b} -I) \hat
        B_{i}^N T_{a \lra
    b}^* +   \hat B_i^N \parens*{T_{a \lra b}^* - I}.\]
    It is easy to verify that non-zero entries of $T_{a \lra b}
    -I$ correspond to the matrix $
    \begin{pmatrix} -1 & 1\\ 1 & -1
    \end{pmatrix}$
    which has a single non-zero eigenvalue of magnitude 2. Thus,
    $\norm*{\Delta_{i,a, b}}_1 \le 4$ where $\norm*{\bullet}_1$ denotes
    the non-normalized trace norm.  Since all the other terms occurring in
    the product $w(\rv{\hat B}^N)$ are random walk matrices (of
    operator norm at most 1), it follows that changing
    each occurrence of $\hat B_i^N$ with $T_{a, b} \hat B_i^N T_{a,b
    }^*$ changes the un-normalized trace at most by 4. Since there
    are $q_i$ occurrences of $\hat B_i^N$ on in $w(\rv {\hat B}^N)$
    replacing all of them, incurs a cost of at most $4q_i$, since we
    consider the normalized trace, we have
    \[ \Abs*{ F(\rv \pi) - F(\rv \pi^{i, a \lra b})} \le \frac{4
    \cdot q_i}{N}.\]
    Thus, we see that we can invoke \cref{lem:perm_conc} with $L_i =
    4q_i/N$. In particular, this implies
    \begin{align*}
        \Pp\sqbr*{ \Abs*{F(\rv \pi) - \Exp F(\rv \pi)} \ge \eta}
        &~\le~2\exp\parens*{- \frac{\eta^2}{2(N-1) \sum_{i=1}^r\frac{
        16 q_i^2}{N^2}}},\\
        &~\le~2\exp\parens*{- \frac{\eta^2 \cdot N}{32 q^2}}
    \end{align*}
    where we used that $q^2 \ge \sum_{i=1}^r q_i^2$. Notice that this
    is equivalent to
    \[ \Pp\sqbr*{ \Abs*{ \tr w(\rv B^N) - \Exp \tr(\rv B^N) } \ge
        \prod_{i=1}^{r} d_i^{q_i} \cdot \eta} \le 2 \exp\parens*{-
        \frac{\eta^2 \cdot N}{32 q^2}}
    \]
    Using \cref{thm:pointwise_convg} and the triangle inequality,
    when $q_i < \girth(G_{N; i})$  for all $i = 1,\ldots, r$ we obtain
    \[ \Pp\sqbr*{ \Abs*{ \tr w(\rv B^N) - \tau(w(\rv t)) } \ge
            \prod_{i=1}^{r} d_i^{q_i} \cdot \parens*{ \frac{3q^2}{N} +
        \eta}} \le 2 \exp\parens*{-
        \frac{\eta^2 \cdot N}{32 q^2}}
    \]
    Notice that when $q_i$ and $d_i$ are fixed we can pick $N$
    sufficiently large to ensure that $q_i$ will satisfy $q_i <
    \girth(G_{N; i})$ for all $i \in [r]$ --- since we have assumed
    $\girth(G_{N, i}) = \omega(1)$.
    Concretely, the choice of $\eta = 8q \sqrt{\frac{\log N}N}$
    ensures that the RHS is bounded by $2/N^2$ while ensuring that
    the deviation measured by the above inequality tends to 0.  This is
    sufficient to show that mixed moments $\tr w(\rv B^N)$ converge
    in probability to $\tau(\rv t)$, concluding our proof.
\end{proof}
\subsection{Spectral Lower Bounds for Permutation Conjugated
    Adjacency
Matrices}\label{sub:spectral_lower_bounds_for_permutation_conjugated_matrices}

First, we note that we can always lower bound the trace $\tr w(\rv
B^N)$ by $\tau(\rv t)$ by a simple covering argument.
\begin{lemma}\label{lem:cover}
    For every word $w = (i_1,\ldots, i_q) \in [r]^q$, we have
    \[ \tr w(\rv B^N) \ge \tau(w(\rv t))\]
    with probability 1 where $\rv B^N$ and $\rv t$ are as defined in
    \eqref{eq:quu_def}.
\end{lemma}
\begin{proof}
    Recall from \cref{sub:pointwise_convergence} that
    $\tau(w(\rv t))$ counts the $w$-colored closed walks emerging
    from the root of $\hat\Tre$, whereas $\tr w(\rv B^N)$ is the
    average number of such walks emerging from a random vertex of
    $\hat G$.

    Let $v$ be any vertex of $\hat G$. We note that as every vertex
    in $\hat G$ has $d_i$ neighbors of color $i$, the connected
    component containing each $v$ can be covered
    by colored infinite tree $\hat \Tre$ with the root mapped to $v$
    --- the covering map here preserves the colors of the edges and
    maps incident edges bijectively.

    This means that every $w$-colored closed walk from the root
    of $\hat\Tre$
    projects to a $w$-colored closed walk from $v$ in $\hat G$.
    Distinct walks have distinct projections as
    every walk in $\hat G$ has a unique lift from the chosen
    root. Consequently,
    \[
        \bigl(w(\rv B^N)\bigr)_{v,v}
        \ge \tau(w(\rv t)).
    \]
    Averaging over $v$ proves the statement.
\end{proof}

With this, we can proceed and prove \cref{thm:general_lowerbound}.
\begin{proof}[Proof of \cref{thm:general_lowerbound}]
    Since $Q$ has non-negative coefficients, so does $Q^* Q$. Thus,
    by invoking \cref{lem:cover} and taking a sum over monomials, we obtain
    \[ \tr\parens*{ \parens*{Q(\rv B^N)^* Q(\rv B^N)}^p} \ge
    \tau\parens*{ \parens*{Q(\rv t)^* Q(\rv t)}^p}.\]
    Writing $\rv d = (d_1,\ldots, d_r)$ for the degree sequence of
    $G_{N; i}$ for $i \in [r]$, we note that the matrix $Q(\rv
    B^N)^* Q(\rv B)$ will always have an eigenvalue $D^2 = Q(\rv
    d)^2$ associated to
    the constant function.
    Thus writing $Q_{\perp 1}(\rv B^N) = \left. Q(\rv
    B^N)\right|_{\perp 1}$, we have
    \[ \Tr\parens*{\parens*{Q(\rv B^N)^* Q(\rv
        B^N) }^p} = D^{2p} + \Tr\parens*{ \parens*{Q_{\perp 1}(\rv
        B^N)^* Q_{\perp 1}(\rv B^N)}^p} \ge N \cdot \tau\parens*{
    \parens*{Q(\rv t)^* Q(\rv t)}^p}.\]
    Equivalently,
    \begin{equation}\label{eq:tr_boud}
        \tr\parens*{ \parens*{Q_{\perp 1}(\rv
        B^N)^* Q_{\perp 1}(\rv B^N)}^p} \ge \tau\parens*{
        \parens*{Q(\rv t)^* Q(\rv t)}^p} - \frac{D^{2p}}N.
    \end{equation}
    We proceed with the choice $p = \max\set*{1, \Big\lceil\sqrt{\log
    N} \Big\rceil}$. Taking $p$-th roots and using the concavity of $x
    \mapsto x^{\frac{1}{2p}}$, we have
    \[ \tr\parens*{ \parens*{Q_{\perp 1}(\rv
        B^N)^* Q_{\perp 1}(\rv B^N)}^p}^{\frac{1}{2p}} \ge \tau\parens*{
        \parens*{Q(\rv t)^* Q(\rv t)}^p}^{\frac{1}{2p}} -
        \frac{D}{N^{\frac{1}{2p}}} =
        \tau\parens*{
        \parens*{Q(\rv t)^* Q(\rv t)}^p}^{\frac{1}{2p}}
        - D \exp\parens*{- \Omega\parens*{\sqrt{\log N}}}
    \]
    For $Q$ a polynomial of fixed degree and $d_i = O(1)$ for all $i
    \in [r]$, we have that $D
    \exp\parens*{-\Omega\parens*{\sqrt{\log N}}} = o_N(1)$. Moreover,
    as $p$ is an
    increasing function of $N$ by
    \cref{lem:norm_def}, we have
    \[ \Abs*{\norm*{Q(\rv
            t)}_{\Cred(\Gam)} -\tau\parens*{
        \parens*{Q(\rv t)^* Q(\rv t)}^p}^{\frac{1}{2p}}} \to
        0\qquad\textrm{as}\qquad N \to \infty.
    \]
    Calling this above correction term $\Delta(N)$ we see that $\Delta(N)
    = o_N(1)$. In particular,
    \[ \tr\parens*{ \parens*{Q_{\perp 1}(\rv
        B^N)^* Q_{\perp 1}(\rv B^N)}^p}^{\frac{1}{2p}}  \ge \norm*{Q(\rv
        t)}_{\Cred(\Gam)} - \Delta(N) - D \exp\parens*{-
            \Omega\parens*{\sqrt{\log
    N}}} = \norm*{Q(\rv t)}_{\Cred(\Gam)} - o_N(1).\]
    Notice now,
    \[ \norm*{Q_{\perp 1}(\rv B^N) } \ge \tr\parens*{
            \parens*{Q_{\perp 1}(\rv
    B^N)^* Q_{\perp 1}(\rv B^N)}^p}^{\frac{1}{2p}} \]
    always. This is because the LHS is the largest singular value of
    $Q_{\perp 1}(\rv B^N)$ whereas the RHS is
    $(\frac{1}{N} \sum_{i = 2}^N \sigma_i^{2p})^{1/(2p)}$ where
    $\sigma_i$ range over the singular values of $Q_{\perp 1}$. In other
    words, the LHS is the $L_\infty$ norm of the vector of
    singular whereas the
    RHS is the $L_{2p}$ norm of the same vector; and we know that
    $L_\infty$ always dominates $L_{2p}$ for $p > 1$. Thus, we conclude
    \[\norm*{Q_{\perp 1}(\rv B^N) } \ge \norm*{Q(\rv
    t)}_{\Cred(\Gam)} - o_N(1).\qedhere\]
\end{proof}
\subsection{Rotating Walks
Reduction}\label{sub:rotating_walks_reduction} 
\begin{proof}[Proof of \cref{prop:rotating_redux}]
    Suppose $\rv U^N$ is a deterministic sequence of permutation matrices.
    \[ \Arot(\rv U_{[2, k-1]}^N) = M^N
        M^{N*}\qquad\textrm{where}\qquad M^N =
    A^N U_2^N A^N U_3^N \cdots U_{k-1}^N A^N.\]
    We set $W_j^N = U^N_1\ldots U_j^N$ for $j \ge 1$.
    Since $W_{j-1}^{N*} W^N_j = U_j^N$ we can write,
    \[ W_1^N A^N W_1^{N^*} \cdots W_{k-1}^N A^N W_{k-1}^{N*} =
        W_1^N A^N U_2^N A^N U_3^N \cdots U_{k-1}^N A^N \cdot
        W_{k-1}^{N*} = W_1^N M^N
    \cdot W_{k-1}^{N^*}.\]
    Consequently,
    \begin{align} (W_1^N A^N W_1^{N^*} \cdots W_{k-1}^N A^N W_{k-1}^{N*})
        (W_1^N A^N W_1^{N^*} \cdots W_{k-1}^N A^N W_{k-1}^{N*})^*
        &~=~ W_1^N
        (M^N M^{N^*}) W_1^{N*},\notag\\
        &~=~ W_1^N \Arot(\rv U^N_{[2, k-1]}) W_1^{N*}.
    \end{align}
    Since $W_1^N = U_1^N$ is a permutation matrix, this means that
    the above matrix has the same spectral law as the rotating walk
    $\Arot(\rv U^N_{[2, k-1]})$.

    We now replace $\rv U^N$ with an iid sequence of uniform
    permutation matrices and keep $\rv W^N$ defined as before. We notice
    that $W_{k-1}^N$ distributes like a
    uniformly random permutation matrix independent from all $\rv
    U^N_1, \ldots, \rv U^N_{k-2}$ thus, replacing $W^N_{k-1}$
    with a uniformly random permutation matrix does not change it's
    spectral law. Repeating the same argument for all $W^N_j$ we see
    that we can assume that each $W_j A^N W_j^{*}$ distributes like
    $B_j^N = U_j^N A^N U_j^{N*}$ whence the statement follows.
\end{proof}


\printbibliography
\end{document}